\documentclass{conm-p-l}

\newcommand{\A}{\mbox {${\mathcal A}$}}
\newcommand{\cP}{\mbox {${\mathcal P}$}}

\newcommand{\CA}{\mbox {${\bf C}{\mathcal A}$}}

\newcommand{\aeta }[1]{\mbox{${}_{#1}\eta $}}
\newcommand{\etaa }[1]{\mbox{$\eta _{#1}$}}

\newcommand{\mapright}[1]{\smash{
   \mathop{\longrightarrow}\limits^{#1}}}

\newcommand{\fig}[1]{\includegraphics{#1}}

\newcommand{\1}{\mbox {${\bf 1}$}}

\newtheorem{ppar}{}[section]

\newtheorem{theo}[ppar]{Theorem}
\newtheorem{prop}[ppar]{Proposition}
\newtheorem{corr}[ppar]{Corolary}
\newtheorem{conj}[ppar]{Conjecture}
\newenvironment{ssec}{\begin{ppar} \rm}{\end{ppar}}

\begin{document}

\title{On Quinn's Invariants of 2-dimensional CW-complexes}

\author{Ivelina Bobtcheva}
\address{Virginia Polytechnic Institute and State University
Blacksburg, VA, 24061}
\email{demeio@anvax1.cineca.it}
\thanks{Partially supported 
by grant from the US National Science Fondation.}

\date{April 1998.}

\begin{abstract}
Given a semisimple stable autonomous tensor category
\A\ over a field $K$,
to any group presentation with finite number of generators we 
associate an element $Q(P)\in K$ invariant under the 
Andrews-Curtis moves. We show that in fact, this is the same
invariant as the one produced by the algorithm introduced by
Frank Quinn in \cite{Q:lectures}. The new definition allows 
us to present a relatively simple proof of the invariance and 
to evaluate $Q(P) $ for some  presentations.
On the basis of some numerical calculations
over different Gelfand-Kazhdan categories,
we make a conjecture  which allows us to
relate the value of $Q(P)$ for two different classes of 
presentations.
\end{abstract}

\maketitle

\section{Introduction}

\label{acmoves}
The original formulation of the Andrews-Curtis conjecture \cite{AC} 
has the following algebraic form. Suppose that
$<x_1,x_2,\dots ,x_n|\; R_1,R_2,\dots ,R_n>$ is a presentation of
the trivial group. Then this can be reduced to the empty presentation
by a finite sequence of operations of the following types:
\begin{itemize}
\item[(i)] The places of $R_1$ and $R_s$ are interchanged.
\item[(ii)] $R_1$ is replaced with $gR_1g^{-1}$, where $g$ is any
           element in the group.
\item[(iii)] $R_1$ is replaced with $R_1^{-1}$.
\item[(iv)] $R_1$ is replaced with $R_1R_2$.
\item[(v)] Adding of an additional generator $y$ and an additional
           relator $y$.
\item[(vi)] The reverse of (v).
\end{itemize}
We will refer to these six operations as AC-moves.
The  Andrews-Curtis conjecture is 
a particular case of the more general deformation conjecture 
stating that a simple homotopy equivalence of 2-complexes can be 
obtained by a deformation through 2-complexes (2-deformation
\footnote{In part of the literature \cite{CM:book} this is actually
called a 3-deformation, since it can be achieved through expansions
and collapses of disks of dimension at most three.}).
It is known that if two complexes $X$ and $Y$ are simple
homotopy equivalent, then for some $k$, there exists a 2-deformation 
from the one point union of
$X$ with $k$ copies of $S^{2}$ to  the one point union of
$Y$ with $k$ copies of $S^{2}$.

 An algebraic formulation of this deformation conjecture comes from
the following result  \cite{Q:handle}: two 2-complexes are 
simple homotopy equivalent if and only if the corresponding group presentations
have the same difference \# relation $-$\# generators, and they can be reduced
to the same group presentation via the moves (i)-(vi) above plus the
additional move
\begin{itemize}
\item[(vii)]
$<x_1,x_2,\dots ,x_n|\; R_1,R_2,\dots ,R_m>$ is replaced with \newline
$<x_1,x_2,\dots ,x_n|\; Q_1,Q_2,\dots ,Q_m>$ if there are words $
w_{i,j}$, $v_{i,j}$ for $1\leq i\leq m$, $1\leq j\leq k(i)$, and 
indices $r(i,j)$, $s(i,j)$ such that for all $i$,
$$
R_{i}(Q_{i})^{-1}=\sqcap 
     _{j=1}^{k(i)}R_{i,j}Q_{i,j}R_{i,j}^{-1}Q_{i,j}^{-1},
$$
\end{itemize}
where $R_{i,j}=w_{i,j}R_{r(i,j)}^{\pm 1}w_{i,j}^{-1}$ and 
$Q_{i,j}= v_{i,j}Q_{s(i,j)}^{\pm 1}v_{i,j}^{-1}$.
Then, in its general form, the Andrews-Curtis  
conjecture states that two presentations which 
have the same difference \# relation $-$ \# generators and are equivalent 
through the moves (i)-(vii), are actually equivalent through the moves
(i)-(vi). 

In \cite{Q:lectures} an algorithm is described  for 
computing topological quantum field
theories on $1+1$-complexes associated to  stable subcategories of 
the  autonomous tensor 
categories 
of Gelfand and Kazhdan  over the finite field $Z_{p}$. 
These two dimensional topological 
theories have the property that the invariant corresponding to $S^2$ is 0,
and therefore give the opportunity of looking for counterexamples of 
the Andrews-Curtis conjecture. The invariants have been approached numerically
 \cite{Q:lectures,thesis,Q:num.pres.}  where 
every CW-complex is being sliced in standard elementary slices. 
Then the invariant of the complex is obtained by composing the morphisms
corresponding to each slice. Up to now no counterexamples have been
detected, and the
analytical value of the invariant even for simple complexes wasn't known.

The goal of this work is to present a ``global'' definition of the
invariant, in sense that, given a group presentation $P$ with a finite number
of generators, we associate to it an element $Q(P)\in K$, which is 
invariant under the AC-moves. There is an explicit formula for $Q(P)$ 
in terms of some kind of  trace (partial trace )
of a given morphism in the underlying semisimple stable autonomous tensor
category \A . 

The paper is structured in the following way. In section 2 
some basic definitions and properties of 
semisimple stable autonomous tensor categories  are presented. In section 3 
the diagrammatic notation , which will be the main tool of
manipulating the expressions is 
described . Section 4 is concerned with the 
definition and properties of partial traces in the category. 
The invariant $Q(P)$ is defined in section 5, and the proof of the
invariance under the AC-moves is presented in section 6. Section 7
reviews the algorithm introduced in \cite{Q:lectures} and  shows that 
it produces the same invariant. In the last section we list some 
conjectures on the dimension functions in the Gelfand-Kazhdan 
categories and on the basis of one of them show that if $P'$ is
a presentation obtained from $P$ by adding a generator $y$ and a 
relator $xyx^{-1}y^{-1}$ (where $x$ is a generator of $P$), and 
if $P''$ is obtained from $P$ by simply
adding a relator $x$, then $Q(P')=NQ(P'')$, where $N$ is the number of 
the simple objects in the category.

I would like to thank Frank Quinn for introducing me into the subject,
and patiently explaining to me many basic concepts and ideas. 
His help and encouragement have been decisive for the completion of
this work.

\section{General category statements}

\begin{ssec}
\A\ is a tensor category if it is supplied with a bifunctor
 $\diamond :\A \times \A \rightarrow \A $ and
 an identity object $\1$, such that for any $A,B,C\in \A$ there are natural 
 isomorphisms
\begin{eqnarray*}
&& \alpha _{A,B,C}:(A\diamond B)\diamond 
 C\rightarrow A\diamond (B\diamond C)\quad  \mbox{(associativity morphisms)}, \\
&&\etaa{A} :\1\diamond A\rightarrow A, \quad\aeta{A} :A\diamond\1
\rightarrow A,\\
&&\gamma =\gamma _{A,B} :A\diamond B\rightarrow B\diamond A
\quad  \mbox{(commutativity morphisms)}.
\end{eqnarray*}
Moreover, $\gamma _{AB}\gamma_{BA} =id_{AB}$. 
These isomorphisms satisfy a number of axioms as listed in \cite{MacL:nat}.

In order to simplify the notation,
in the future we will often write  $A B$ instead of $A\diamond B$.
\end{ssec}

\begin{ssec} 
\label{bracketing}
Let \cP\ be the free $(\otimes ,I)$-algebra on a single symbol $O$.
The number of $O$'s in an element $T$ in \cP\ will be called
its length. 
Let now \A\ be a tensor category,. Following \cite{KL}, 
we define a category $\cP\circ \A$, with objects 
$T[\underline{B}]\equiv T[B_{1}, B_{2},\dots ,B_{n}]$, where $T$ 
is an object of length $n$ 
in \cP, and $\underline{ B}=(B_{1}, B_{2},\dots ,B_{n})$ is
a sequence of objects in \A. The morphisms have the form
$\pi [f_1, f_2, \dots , f_n]:T[B_{1}, B_{2},\dots ,B_{n}]\rightarrow
S[A_{1}, A_{2},\dots ,A_{n}]$, where $\pi $ is a bijection from the
$O$'s in $T$ to the $O$'s in $S$, i.e. an element in
the symmetric group ${\bf S}_n$, and
$f_i:B_{\pi ^{-1}(i)}\rightarrow A_i$ are maps in \A . We use the 
following convention for the special maps:
\begin{eqnarray*}
&&T[f_1, f_2, \dots , f_n]\equiv \iota (n)[f_1, f_2, \dots , f_n]:
   T[B_{1}, B_{2},\dots ,B_{n}]\rightarrow
      T[A_{1}, A_{2},\dots ,A_{n}],\;\\
&&\pi\equiv \pi [id_{A_1},id_{A_2},\dots ,id_{A_n}]:
    T[A_{\pi(1)}, A_{\pi(2)},\dots ,A_{\pi(n)}]\rightarrow
      S[A_{1}, A_{2},\dots ,A_{n}],
\end{eqnarray*}
where $f_i:B_i\rightarrow A_i$, and  $\iota (n)$ denotes the identity element 
${\bf S_n}$.
We will refer to the objects in $\cP\circ \A$ as bracketings
in \A, and to the maps of type 
$\pi=\pi [id_{A_1},id_{A_2},\dots ,id_{A_n}]$ as permutations. 
\end{ssec}

\begin{ssec}
\label{tenstrPB}
$\cP\circ \A$ possesses a structure of a tensor category
with identity object $I[\; ]$, and product $\otimes $ defined in the
following way:
$$
 S[B_{1}, B_{2},\dots ,B_{n}] \otimes  T[B_{n+1}, B_{2},\dots ,B_{m}]=
   (S\otimes T)(B_{1}, B_{2},\dots ,B_{n+m}).
$$
The associativity and commutativity morphisms are defined to be the
corresponding permutations. 
 Then the coherence theorem of MacLane \cite{MacL:nat}
is equivalent to the assertion that there is a strict tensor functor
$\cP\circ \A\rightarrow \A$. We write 
$T(B_{1}, B_{2},\dots ,B_{n})$ and 
$\pi (f_1, f_2, \dots , f_n)$ for the images of $T[B_{1}, B_{2},\dots ,B_{n}]$
and $\pi [f_1, f_2, \dots , f_n]$ under this functor; for example 
$((O\otimes O)\otimes I)(a,b)=(a\diamond b)\diamond \1$. The images of
the permutations will be called again permutations, and  
we will denote them with the corresponding 
cyclic decomposition. For 
example, $(1,2)(3,5,4): T(A,B,C,D,G)\rightarrow S(B,A,D,G,C)$.

In the case when $B_{i}=B$ for any $i$, we will use the convention
$T(B^k)=T(B_{1}, B_{2},\dots ,B_{k})$. Then the permutations
 define an action of ${\bf S}_{k}$ on $T(B^k)$, and 
we denote the corresponding representation of ${\bf S}_{k}$
 as $\rho [T(B^k)]$.
\end{ssec}

\begin{ssec}
Let $P\A $ be the category with the same objects as $\cP\circ \A$,
but has different set of morphisms: for any morphism 
$T(\underline{B})
\rightarrow R(\underline{ B})$ in \A\, there is a 
 morphism $T[\underline{ B}]
\rightarrow R[\underline{ C}]$, which will be denoted in the same way.
The composition is obvious. Then there are functors $ \cP\circ \A
\rightarrow P\A$ and
$P\A\rightarrow \A$, whose compositions is exactly the functor 
$ \cP\circ \A\rightarrow \A$ given by the coherence result of MacLane.
\end{ssec}

\begin{ssec}
\label{defn.satc}
An {\em autonomous tensor} category in the terminology of \cite{Shum} is a
tensor category \A\ in which to every object $A$ is  assigned an object
$A^*\in \A$ and morphisms $\Lambda _{A}:\1\rightarrow A^{*}A$ (coform) and
$\lambda _{A}:AA^{*}\rightarrow \1$ (form) such that the compositions
\begin{eqnarray*}
&&A^*  \mapright{\etaa{A^*} ^{-1}} \1\, A^*  \mapright{\Lambda _{A}\diamond 
id_{A^*}} 
  (A^{*}A)A^*  \mapright{\alpha } A^*(AA^*)  
  \mapright{id_{A^*}\diamond \lambda _{A}} A^*\,\1  \mapright{\etaa{A^*}} 
  A^*, \\
&&A  \mapright{\aeta{A} ^{-1}} A\, \1  \mapright{id_{A}\diamond 
\Lambda _{A}} 
  A(A^{*}A)  \mapright{\alpha ^{-1}} (AA^{*})A  
  \mapright{ \lambda _{A}\diamond id_{A}}  \1 \,A 
  \mapright{\etaa{A}} A   
\end{eqnarray*}
act as identities.
In any such category there is a canonical isomorphism
$w_{A}:A^{**}\rightarrow A$, for which
$
\gamma_{A^{*}A}\circ \Lambda _{A}=(w_{A}\diamond id)\circ\Lambda _{A^{*}}
$
We call an autonomous tensor category \A\ {\em stable} 
if $(A^{*})^{*}=A$  and $w_{A}=id_{A}$.
\end{ssec}

\begin{ssec}
\label{defn.semis}
A category \A\ is called {\em semisimple} if it is abelian over a 
field $K$, and if there is a 
finite subset of objects $\Sigma $ in \A, such that every other 
objects is isomorphic to a direct sum of objects in $\Sigma $, and
for any $a,b\in\Sigma$, 
$hom _{\A}(a,b)\simeq \left\{ \begin{array}{c} 0,\mbox{ if } a\neq b,\\
                                               K,\mbox{ if } a= b.
                               \end{array}
                                 \right. $.
In general, we will use small latin
letters $a,b,c,\dots $ to denote objects in $\Sigma $, and capital
letters to indicate an arbitrary object in the category. We also
use the notation $F(A,B)=hom_{\A}(A,B)$. 

The fact that an object $A$ is isomorphic to a direct sum of objects 
in $\Sigma $ means that for any $a \in \Sigma$ there exists a basis
$\{\epsilon _i(a,A):a\rightarrow A\}_{i}$ of $F(a,A)$ and a basis
$\{\epsilon _i(a,A)^*:A\rightarrow a\}_{i}$ of $F(A,a)$, which
are dual in sense that 
$\epsilon _i(a,A)^*\circ \epsilon _j(a,A)=\delta _{i,j}id_a$. Moreover, 
$\sum _{a,j}\epsilon _j (a,A) \circ \epsilon _j(a,A)^*= id _{A}$. 
Then,
given the objects $A,B$ and $C$, the map
\begin{eqnarray*}
   \nabla : &\oplus _{b\in\Sigma} F(a,A\diamond (b\diamond C))\otimes
                        F(b,B)\rightarrow &F(a,A\diamond (B\diamond C)),\\
      &\psi \otimes \phi\rightarrow &(id_A\diamond (\phi\diamond id_C))
        \circ\psi
\end{eqnarray*}
is an isomorphism. In fact, the inverse of $\nabla$ is given by
the map
\begin{eqnarray*}
&\nabla^{-1}:F(a,A\diamond (B\diamond C))\rightarrow 
     &\oplus _{b\in\Sigma} F(a,A\diamond (b\diamond C))\otimes F(b,B),\\
&\phi \rightarrow &\oplus _{b\in\Sigma ,i}
   [(id_A\diamond (\epsilon _i(b,B)^*\diamond id_C))\circ \phi]\otimes 
    \epsilon _i(b,B).
\end{eqnarray*}

In a semisimple tensor category for any simple object $b$, 
the representation $\rho [T(a^k)]$ 
from \ref{tenstrPB}, induces representations 
$\rho _{b}[T(a^k)]$ of ${\bf S_{k}}$ on the space $F(b, T(a^k))$.
\end{ssec}

\begin{ssec}
\label{perp.cat.} 
The concrete  categories, used 
in the numerical computations of  quantum invariants of 
2-dimensional CW-complexes, are the ones studied in \cite{GK}.
These are semisimple tensor  categories, 
defined over the finite field $K=Z_p$, and are constructed as a quotient
of the category of some modular finite-dimensional
modules of a semisimple algebraic group $G$, in the case when
$ p$ is greater
then the Coxeter number of the corresponding Lie-algebra. 
If \A\ is such Gelfand-Kazhdan category,
its simple objects
correspond to the simple modules with highest weights inside the 
fundamental alcove. We remind that {\em involution} on a set $S$ is a 
bijection of $S$ into itself whose square is equal to the identity.
Then on the set of simple objects $\Sigma $ 
there is an involution denoted by $a\rightarrow a^{*}$ , 
such that $a^{*}$ is isomorphic to the dual module of $a$. A form and
coform can be introduced, and this makes \A\ an autonomous tensor 
category. But this category is not always stable, i.e. the condition
$
 \gamma _{a,a^{*}}\circ \Lambda _{a^{*}}  =\Lambda _{a}:\1
\rightarrow a^{*}\diamond a.
$
is not always satisfied. The simplest examples come from 
the categories
of modular representations of $sl(2)$. In this case if $a$ is an odd
highest weight module in $\Sigma $,
$
 \gamma _{a,a^{*}}\circ \Lambda _{a^{*}}=-\Lambda _{a}.
$
But the subcategory generated by the even highest weights 
modules is stable. 

The general situation is not very different from the one of $sl(2)$, 
and is described, for example, in \cite{thesis}. It is shown
that there is always a subcategory of 
\A\ which is stable, for example the one generated by the
modules with highest weights in the root lattice. In general, 
there are number of 
sublattices of the weight lattice which produce subcategories with the desired
property. 
\end{ssec}

\section{Diagrammatic notation}

The diagrammatic notation used here is based on the coherence results in 
\cite {KL,Shum,FY}. The particular case of autonomous tensor 
categories is studied in \cite{KL}, and the results there are extended 
to tortile categories in \cite {Shum,FY}. The version we use here 
is the one described in \cite{Y},
specified to the case of stable autonomous tensor category.

The idea behind the diagrammatic notation 
is that to each composition of $\A$-morphisms there corresponds
a labeled link diagram with coupons. Then there are rules for 
transforming these diagrams, and the statement is that if by such 
transformations the diagrams corresponding to $\psi 
_{1}:A\rightarrow B$ and $\psi 
_{2}:A\rightarrow B$ can be reduced to the same diagram,
then $\psi _{1}=\psi _{2}$. In this sense, a morphism in \A\
is determined by its diagram.

\begin{prop} 
\label{invprop}
Suppose that the group $Z_2\times Z_2=<\theta _{1}, \theta _{-1}|\;
\theta _{1}^{2}, \theta _{-1}^{2}>$ acts on a set $S$. Let 
$S_{i}$, $i=\pm 1$, be the set of fixed points of the subgroup generated by 
$\theta _{i}$, and suppose that $S_{1}\cap S_{-1}=\emptyset $. Then if
$\tilde{S}=S\setminus (S_{1}\cup S_{-1})$ and ${\mathcal O}$
is any orbit, we have that either ${\mathcal O}\setminus \tilde{S}$ is
a two elements set, or ${\mathcal O}\subset \tilde{S} $. 
\end{prop}
The proof is straightforward.

\begin{ssec}
The {\em coupon
category} \CA\ associated to a stable autonomous tensor category 
$(\A , \alpha ,\gamma ,\eta ,\lambda ,\Lambda ,*)$
is defined 
as a quotient of another category $\tilde {\CA }$.
The objects of $\tilde {\CA }$ are finite sequences of objects in \A .
Below, we will use the notation $\underline{A}$ both for the sequence 
$(A_{1},A_{2}, \ldots , A_{l})$ of length
$|\underline{A}|=l$, and for the set of objects $\{A_i\}_{i=1}^{l}$ in \A\ .
In particular, $\underline{A}\sqcup \underline{B}$ denotes the disjoint
union of the two sets, and $\underline{A}\setminus \underline{B}$
denotes the complement of $\underline{B}$ in $\underline{A}$.

We introduce a 
$*$ operation on the objects in $\tilde {\CA }$: $
\underline{A}=(A_{1},A_{2}, \ldots , 
A_{l})\rightarrow \underline{A}^{*}=(A_{l}^{*},A_{l-1}^{*}, \ldots , A_{1}^{*})$.
Then a morphism $(F,\theta ):\underline{X}\rightarrow \underline{Y}$ in 
$\tilde {\CA }$ is given by a 
finite set $F=\{f^i:T^{i}[\underline{A}^{i}]\rightarrow
 R^{i}[\underline{B}^{i}]\}_{i}$ of morphisms 
in  $P\A $ , and an involution $\theta $ on the set  
$S_{F}=\underline{X}\sqcup \underline{Y}^{*}\sqcup 
(\sqcup _{i}(\underline{A}^{i}{}^{*}\sqcup \underline{B}^{i}))$ such that 
$\theta (C)^{*}=C$ and $\theta$ is free of fixed points. The elements of
$\underline{A}^{i}$ and $\underline{B}^{i}$ are called correspondingly the 
lower and the upper ends of $f^{i}$. 
The identity morphisms in $\tilde {\CA }$, 
$id :S\rightarrow S$ are given by $(\{\}, *)$. 

Let now
$(F,\theta ):\underline{X}\rightarrow \underline{Y}$ and 
$(G,\eta ):\underline{Y}\rightarrow 
\underline{Z}$  be two morphisms in $\tilde {\CA }$. To define their composition,
we would like to take $F\sqcup G$ with the involution map obtained by
identifying the two copies of $\underline{Y}$, but
in this way one obtains an involution, which may contain fixed points.
To avoid this, identity morphisms are introduced as follows. 

Let $S'= S_{F}\setminus \underline{Y}^{*}$
and $S''=S_{G}^{*}\setminus \underline{Y}^{*}$. We
 extend $\theta$ to an involution on $S=
S'\sqcup \underline{Y}^{*}\sqcup S''$ by defining $\theta$ to be the 
identity on the elements in $S''$. In a similar way we define an 
involution $\eta '$ on $S$ as $\eta' (x)=\left\{ 
\begin{array}{c} \eta (x^*)^{*} \mbox{ if } x\in \underline{Y}^{*}\sqcup S''\\
                   x \mbox{ otherwise }
 \end{array}\right. $. Then $\theta$ and $\eta '$ define an action of
 $Z_{2}\times Z_{2}$ on $S$, which satisfy the requirements of
proposition \ref{invprop}, where $\tilde{S}=\underline{Y}^{*}$, $S_1=S'$,
and $S_{-1}=S''$. Then, according to the proposition, the orbits under this
action are divided in two - the ones which are entirely contained in 
$\underline{Y}^{*}$, and ones which have exactly two elements
outside this set. 
Let ${\bf E}$ be the subset of orbits 
 which are entirely contained in $\underline{Y}^{*}$. 
If ${\mathcal O}\in {\bf E}$, let $A_{\mathcal O}$ 
be the smallest element in the
orbit ( with respect to the ordering which comes with $\underline{Y}^{*}$), and
let ${\bf T}=\{id_{A_{\mathcal O}}\}_{{\mathcal O}\in {\bf E}}$. 
Then composition of $(F,\theta )$ and $(G,\eta )$ is given by
$(F\sqcup G\sqcup T, \xi )$, where $\xi$ is an involution on the set
$S'\sqcup (S'')^{*}\sqcup 
(\sqcup _{{\mathcal O}\in {\bf E}}\{A_{\mathcal O}^{*},A_{\mathcal O}
\})$ defined in the following way. Let $\imath :S'\sqcup (S'')^{*}
\stackrel{id\sqcup *}{\longrightarrow}
 S'\sqcup S''\subset S$. Then if 
$A$ is an element in $S'\sqcup (S'')^{*}$, we define $\xi (A)=\imath ^{-1}(B)$, 
where $B$ is the unique element in $S'\sqcup S''$  which belongs
to the same $Z_{2}\times Z_{2}$-orbit as $\imath (A)$.
Moreover, $\xi (A_{\mathcal O})=A_{\mathcal O}^{*}$. 
\end{ssec}

\begin{ssec}
\label{carel}
The coupon category \CA\ is defined to be the quotient of $\tilde{\CA}$ under the 
following relations.
Let $(F,\theta )$ be a morphism in $\tilde{\CA}$. Then the following
replacements between morphisms in $F$ can be made:
\begin{itemize}
\item[(a)]  $id _{A}:O[A]\rightarrow O[A] $, whose ends 
   are not mapped into each other, can be removed. Then $\theta $ changes to
   $\theta '$ where $\theta '(X)$ is equal to 
   $\theta (X)$ if $X$ is not mapped into some of the ends of $id _{A}$, 
   otherwise $\theta '(\theta (A))=\theta (A^{*})$.
\item[(b)] If $f:T[\underline{A}]\rightarrow Q[\underline{B}]$ 
   and $g:Q[\underline{C}]\rightarrow R[\underline{D}]$ are 
   such that $|\underline{B}|=|\underline{C}|$ and the
    upper ends of $f$ are mapped into the lower 
    ends of $g$, i.e. $\theta (B_{i})=C_{i}^{*}$, then $f,g$ can be replaced
    with $f\circ g$. 
\item[(c)] $f:T[\underline{A}]\rightarrow Q[\underline{B} 
    ]\;\longleftrightarrow \;f'=\iota (l)^{-1}\circ f\circ \iota (k):
    T'[\underline{A}]\rightarrow Q'[\underline{B}]$, where
    $k=|\underline{A}|$, $l=|\underline{B}|$, and $\iota (k):
    T'(\underline{A})\rightarrow T(\underline{A})$ and $\iota (l):
    Q'(\underline{B})\rightarrow Q(\underline{B})$ are the identity
    permutations as defined in 
    \ref{bracketing}.
\item[(d)]  $\gamma _{A,B}:(O\otimes O)[A,B]\rightarrow (O\otimes 
   O)[B,A]\;\longleftrightarrow \;
   id :(O\otimes O)[A,B]\rightarrow (O\otimes O)[A,B]$.
\item[(e)] $\{f:T[\underline{A}]\rightarrow T'[\underline{B}],\, 
    g:Q[\underline{C}]\rightarrow Q'[\underline{D}]\}
    \;\longleftrightarrow \;
    g\diamond f: T[\underline{A}]\otimes Q[\underline{B}]\rightarrow 
    T'[\underline{C}]\otimes Q'[\underline{D}] $.
\item[(f)] $\eta _{A}:(O\otimes I)[A] \rightarrow O[A] 
    \;\longleftrightarrow \;
    id_{A}: O[A]\rightarrow O[A]$, and the analogous statement for 
   ${}_{A}\eta $.
\item[(g)] $\lambda _{A}:(O\otimes O)[A,A^{*}] \rightarrow 
   I[\{ \}]\;\longleftrightarrow \;
   id_{A^{*}}:O[A^{*}]\rightarrow O[A^{*}]$.
\item[(h)]$\Lambda _{A}:I[\{ \}]\rightarrow(O\otimes O)[A^{*},A] 
   \;\longleftrightarrow \;
   id_{A}:  O[A]\rightarrow O[A]$.
\end{itemize}
\end{ssec}   
   
\begin{ssec} 
The category \CA\ is a stable autonomous tensor category.
The product $\underline{A}\otimes \underline{B}$ is simply the sequence 
obtained by putting the elements of $\underline{B}$ after the elements of 
$\underline{A}$,
the identity object is the empty set, and the associativities, 
$\eta _{\underline{A}}$ and ${}_{\underline{A}}\eta $ are the identity morphisms. 
$\gamma $, $\lambda $ and $\Lambda $ are given by $(\{ \}, *)$.
Moreover, 
there is a strict tensor 
functor $P\A\rightarrow \CA$ which maps $T[\underline{ B}]$ into $
\underline{ B}$ and 
$f:T[\underline{ B}]\rightarrow R[\underline{ C} ]$ into 
$(\{f\},*):\underline{ B}\rightarrow \underline{ C}$. 
Then the coherence result which underlines the use of the diagrams is that 
if $\psi ,\phi:
T(\underline{ B})\rightarrow R(\underline{ C})$ are two morphisms in \A\ 
which have the same image in \CA, then $\psi =\phi$. In this way, given 
$T,R\in \cP $, every morphism $(F,\theta ):\underline{ B}\rightarrow 
\underline{ C}$ determines a unique morphism 
$T(\underline{ B})\rightarrow R(\underline{ C})$ in \A, which also will be
denoted as $(F,\theta )$.
\end{ssec}   

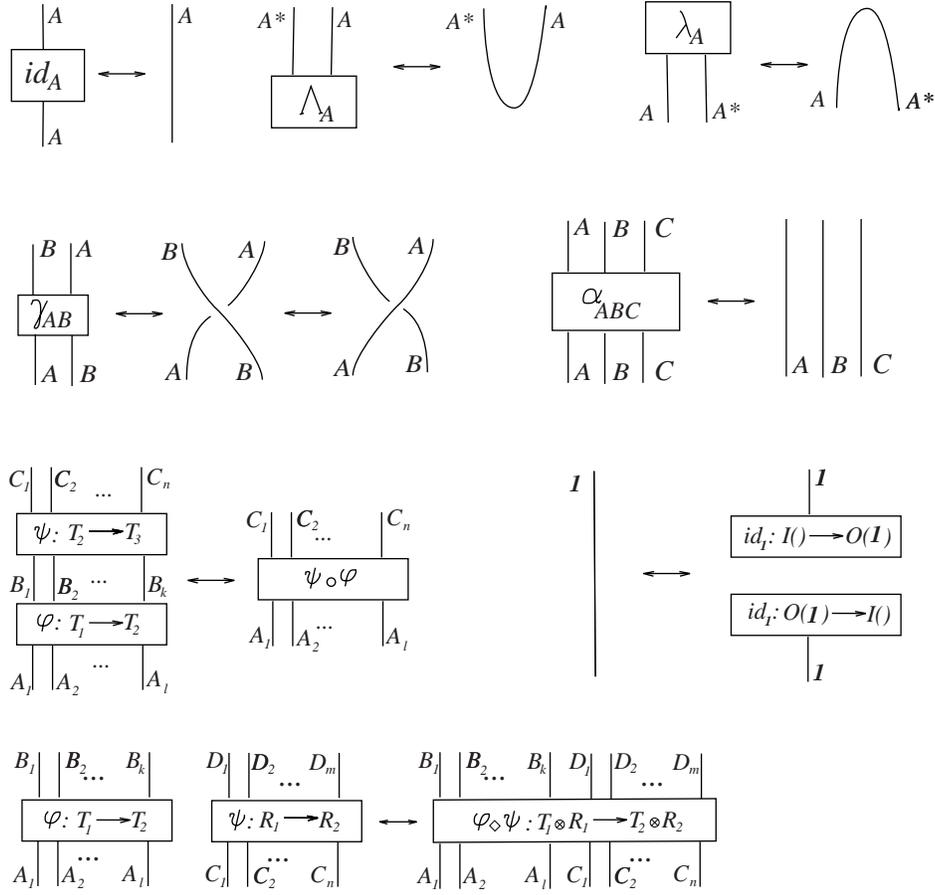
\begin{figure}[h]
\setlength{\unitlength}{1cm}
\begin{center}
\begin{picture}(15,12)
\put(0,0){\fig{relCA.pict}}
\end{picture}
\end{center}
\caption{Diagram presentation of the relation \ref{carel}}
\label{relations}
\end{figure}

\begin{ssec}
A morphism  $(F,\theta ):\underline{X}\rightarrow \underline{Y}$ 
in \CA\  
can be described by a diagram in the following way. 
First we draw a set of $|\underline{X}|+ |\underline{Y}|$ points labeled by the
corresponding elements in $\underline{X}$ 
and $\underline{Y}$, and for each morphism 
$f:T[\underline{A}]\rightarrow
R[\underline{B}]$ we draw a rectangle (coupon) 
labeled by $f$ with $|\underline{A}|$ lines 
labeled by $A_i$'s,
attached to its lower side, and $|\underline{B}|$ lines labeled by 
$B_j$'s, 
attached 
to its upper side. Then we connect with lines every two ends of coupons 
and labeled points which are mapped into each other by the 
involution $\theta $. In this way one looses any 
information about the elements in \cP\ associated with the 
morphism of a coupon except their length, 
but note that according to \ref{carel} (c) this is the only essential 
one. An example for the type of diagram  obtained in this way is 
presented in figure \ref{gendiag}.  Then the equivalence relations in
\ref{carel} imply that one can perform the local diagram moves represented on 
figure \ref{relations}. Note that if in the diagram there is a line which 
starts and finishes in a coupon and is labeled by the
identity object \1, this line can be removed.

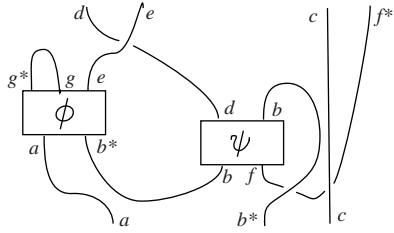
\begin{figure}[h]
\setlength{\unitlength}{1cm}
\begin{center}
\begin{picture}(5,3)
\put(0,-0.3){\fig{fig1.pict}}
\end{picture}
\end{center}
\caption{General diagram.}
\label{gendiag}
\end{figure}
\end{ssec}

\begin{figure}[h]
\setlength{\unitlength}{1cm}
\begin{center}
\begin{picture}(12,3)
\put(0,-0.3){\fig{fig7.pict}}
\end{picture}
\end{center}
\caption{Diagram presentation of the defining identities for the dual 
          basis.}
\label{id.basis}
\end{figure}
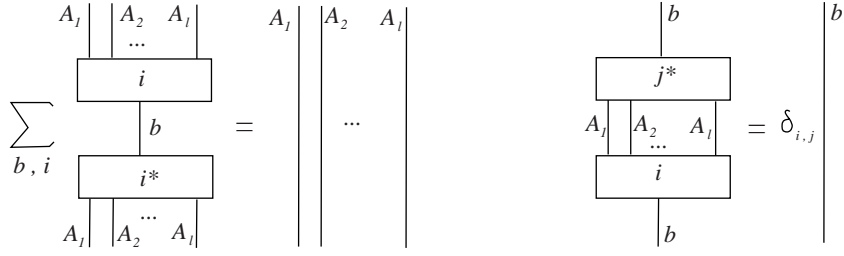

\begin{ssec}
We introduce a special notation for the basic elements of the homomorphism 
spaces.
If $\{\epsilon _{i}\}_{i}$ is a basis of $F(b, 
T(A_{1},A_{2},\dots , A_{l}))$ and  $\{\epsilon _{i}^{*}\}_{i}$ is its
dual basis of $F(T(A_{1},A_{2},\dots , A_{l}),b)$ (see \ref{defn.semis}), 
we represent these basic elements in the
diagrams as it is shown in figure \ref{basis} (a) and (b). In the special and 
most often used case of a basis for $F(a,bc)$, where $a,b,c\in 
\Sigma$, the notation will be simplified by replacing the coupon with a
vertex labeled by $i$.  Figure
\ref{id.basis} represents the defining identities of the dual basis.
From the discussion in \ref{defn.semis} it follows 
that if $\{\zeta _k\}_{k}$ is a basis 
for $F(a,bC)$ and $\{\epsilon _i\}_{i}$ is a basis 
for $F(b,AB)$, then as one varies $b,i,k$, the morphisms
$(\epsilon _{i}\diamond id_{C})\circ \zeta _{k}$ form a basis of
$F(a,(AB)C)$ (figure \ref{basis} (c) ). 
\end{ssec}
\begin{figure}[h]
\setlength{\unitlength}{1cm}
\begin{center}
\begin{picture}(12,2)
\put(0,-0.3){\fig{fig6.pict}}
\end{picture}
\end{center}
\caption{Diagram presentation of basic elements.}
\label{basis}
\end{figure}
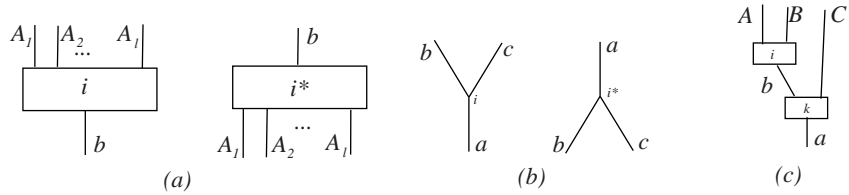

\begin{ssec}
\label{conj.basis}
 Given a morphism
$\varphi$ of $F(A,B)$, we define
$$
\hat{\varphi}= \eta _{A^{*}}(\lambda _{B^{*}}\diamond id_{A^{*}})
\alpha _{B^{*},B,A^{*}}^{-1}
(id_{B^{*}}\diamond (\varphi \diamond id_{A^{*}}))
(id_{B^{*}}\diamond\Lambda _{A^{*}})({}_{B^{*}}\eta )^{-1}\in F(B^*,A^*).
$$
The corresponding diagram is presented in figure \ref{hatbasis}. 
It is easy to see that if $\{\epsilon _i\}_i$ is a basis for
$F(A,B)$, then $\{\hat{\epsilon }_{i}\}_i$ is a basis for $F(B^*,A^*)$ with dual
 $\{\hat{\epsilon ^{*}_{i}}\}_i$.

\begin{figure}[h]
\setlength{\unitlength}{1cm}
\begin{center}
\begin{picture}(4,2)
\put(0,-0.3){\fig{fighatbasis.pict}}
\end{picture}
\end{center}
\caption{Diagram presentation of $\hat{\varphi}$.}
\label{hatbasis}
\end{figure}
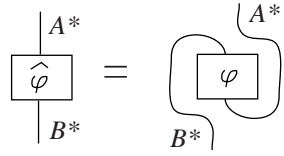
\end {ssec}

\begin{figure}[h]
\setlength{\unitlength}{1cm}
\begin{center}
\begin{picture}(9,1.5)
\put(0,-0.3){\fig{figrank.pict}}
\end{picture}
\end{center}
\caption{Rank of a morphism and rank of an object.}
\label{fig.rank}
\end{figure}
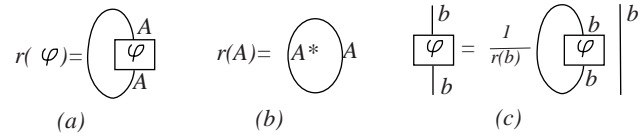

\begin{ssec}
\label{rank}
If $\varphi :A\rightarrow A$ is a morphism in \A , define the rank
of $\varphi$ to be the morphism
$r(\varphi)=\lambda _{A^*}\circ (id_{A^{*}}\diamond \varphi )\circ 
\Lambda _{A}$, and
$r(A)\equiv r(id_{A})$. We think of $r(\varphi)$ as an 
element in the field via the standard isomorphism $F(\1 ,\1)\sim 
K$. The diagrammatic notation for the rank of a 
morphism  
is presented in figure \ref{fig.rank} (a). As an easy application of the 
diagrams one can see that
\begin{itemize}
\item[(a)] $r(A)=r(A^{*})$ ;
\item[(b)] $r(T(A_{1},A_{2}, \ldots , 
       A_{l}))=r(A_{1})r(A_{2}) \ldots r(A_{l})$;
\end{itemize}
(a) allows as to use the simplified notation on 
figure \ref{fig.rank} (b) for a rank of an object.
The identity in figure  \ref{fig.rank} (c) follows from the fact that 
for any simple object $b$, $F(b ,b)\sim K$. 
\end{ssec}

Let $\{\epsilon _{i}(a,bc)\}_{i}$ be a basis for $F(a,bc)$, and 
$\{\epsilon _{i}(a,bc)^*\}_{i}$ 
be the corresponding dual basis. Then the morphisms $\bar{\zeta 
}_{i}(c,b^{*}a)=
\etaa{c}(\lambda _{b^*}\diamond id_{c})(id_{b^*}\diamond \epsilon _{i})$
and $\zeta _{i}(c,b^{*}a)=(id _{b^*}\diamond \epsilon _{i}^*)
(\Lambda _{b}\diamond 
id_{c})\etaa{c}^{-1}$ represent sets of elements in $F(b^{*}a,c)$ and
$F(c,b^{*}a)$ correspondingly. 

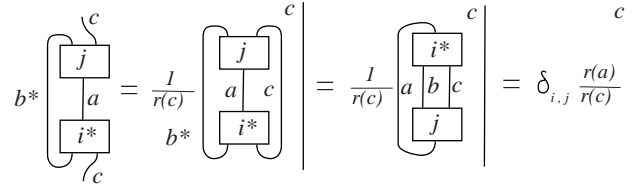
\begin{figure}[h]
\setlength{\unitlength}{1cm}
\begin{center}
\begin{picture}(9,2.2)
\put(0,-0.3){\fig{figdbasis.pict}}
\end{picture}
\end{center}
\caption{Proof of proposition \ref{dbasis}.}
\label{fig.dbasis}
\end{figure}

\begin{prop} \label{dbasis} 
$\{\zeta_{i}(c,b^{*}a)\}_i$ forms a basis for $F(c,b^*a)$ with a dual basis
given by 
$\zeta _{i}(c,b^{*}a)^*=\frac{r(c)}{r(a)}\bar{\zeta}_{i}(c,b^{*}a)$. 
In particular,
$\sum _{c, i} \frac{r(c)}{r(a)}(\zeta _{i}(c,b^{*}a)
\circ\bar{\zeta}_{i}(c,b^{*}a))=id_{b^*a}$.
\end{prop}
The proof is presented in figure \ref{fig.dbasis}.

\begin{ssec}
To make it easier in the future to construct diagrams and to study
parts of them, we introduce the notion of an {\em open morphism}, 
correspondingly {\em open diagrams}. An open morphism in \CA\ is 
a morphism $(F,\theta ):\underline{A}=(A_{1},A_{2}, \ldots ,A_{l})\rightarrow 
\underline{B}=(B_{1},B_{2}, \ldots ,B_{k})$ where some subset 
$In(F)$ of the lower 
ends is labeled as $in$-ends, and some subset $Out(F)$ of the upper 
ends is labeled as $out$-ends. We refer to the set of $in$- and 
$out$- ends as {\em free} ends. 

If 
$(F,\theta):\underline{A}\rightarrow \underline{B}$ and 
$(G,\xi ):\underline{C}\rightarrow \underline{D}$ are two open 
morphisms, and $S\subset Out(F)\cap In(G)$, we define the product of 
the two morphisms over $S$ to be the open morphism obtained by mapping 
the $out$-ends of $(F,\theta)$ labeled by $S$ into the $in$-ends of 
$(G,\xi )$ labeled by $S$ as shown on figure \ref{fig.prodopen} (a). 
To be more precise,
let $\underline{B}'=
\underline{B}\setminus S$ and $\underline{C}'=\underline{C}\setminus 
S$. Then
$(F,\theta)\circ _{S}(G,\xi ):\underline{A}\otimes 
\underline{C}'\rightarrow \underline{B}'\otimes 
\underline{D}$ with $In( G\circ _{S}F)=In(F)\sqcup 
(In(G)\setminus S)$ and
$Out( G\circ _{S}F)=(Out(F)\setminus S)\sqcup Out(G)$ 
is the following composition in \CA :
$$
(F,\theta)\circ _{S}(G,\xi )=(id_{\underline{B}'}\otimes  (G,\xi ))
\circ (\{\},\epsilon )\circ 
((F,\theta)\otimes id_{\underline{C}'}), 
$$
where  $\epsilon :\underline{B}\sqcup \underline{C}'\rightarrow
(\underline{B}'\sqcup \underline{C})^{*}$
is the involution map sending $S\subset \underline{B}\rightarrow 
(S\subset \underline{C})^{*}$.

\begin{figure}[h]
\setlength{\unitlength}{1cm}
\begin{center}
\begin{picture}(16,6.5)
\put(0,-0.3){\fig{figopen.pict}}
\end{picture}
\end{center}
\caption{Product, closure and reverse of open morphisms.}
\label{fig.prodopen}
\end{figure}
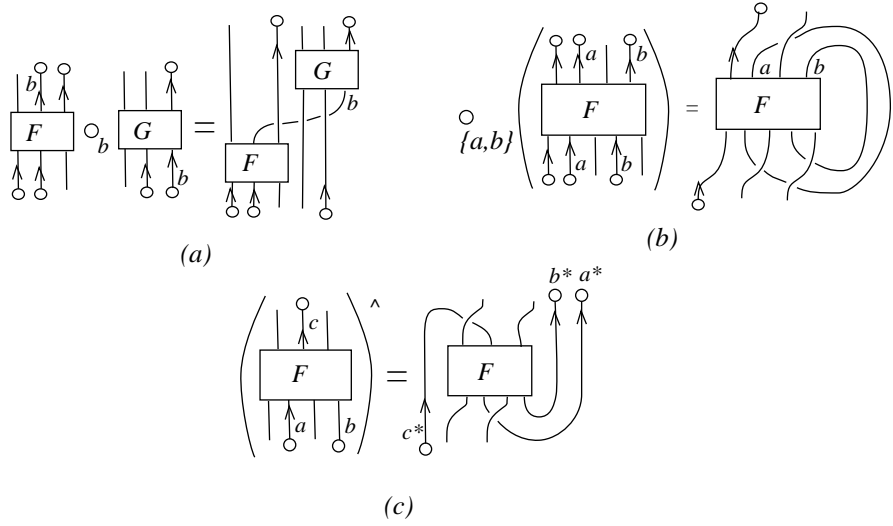

In the diagrams we indicate the $in$-ends with a circle and an arrow 
which comes out of it, and the $out$-ends with a circle and an arrow 
which goes into it. We 
will think of an open diagram as a sub-diagram of a bigger one, 
obtained by cutting out in the places of the free ends. For this 
reason, often the free ends will be left hanging  without 
putting them in a specific place 
among the lower or the upper ends of the diagram, 
indicating in this way that the only important thing is how these will
glue with the free ends of other open diagrams in the final ``closed'' one. 
\end{ssec}

\begin{ssec}
Given an open morphism 
$(F,\theta):\underline{A}\rightarrow \underline{B}$
 we define the closure and the reverse
of $(F,\theta)$. Given 
$S\subset In(F)\cap Out(F)$, the closure
of $(F,\theta)$ over $S$ is defined  to be the open morphisms
$\circ _{S}(F,\theta )=(F,\theta ')$, where $\theta '$ is the extension of
$\theta$ obtained by mapping into each other $S\subset \underline{A}$ and
$S^*\subset \underline{B}^*$ as shown in figure \ref{fig.prodopen} (b).
By definition $In(\circ _{S}(F))=In(F)\setminus S$ and
$Out(\circ _{S}(F))=Out(F)\setminus S$.

The reverse of  $(F,\theta)$ is defined to be a morphism
$(F,\theta)^{\hat{}}=(\hat{F},\hat{\theta})
:Out(F)^{*}\otimes (\underline{A}\setminus 
In(F))\rightarrow (\underline{B}\setminus 
Out(F))\otimes In(F)^{*}$ with $In(\hat{F})=Out(F)^{*}$ and
$Out(\hat{F})=In(F)^{*}$. $(F,\theta)^{\hat{}}$ is obtained from 
$(F,\theta )$ by exchanging the
places of the $in$- and $out$- ends through composition with the corresponding
form and coform in \CA\ as it is shown in figure 
\ref{fig.prodopen} (c).
\end{ssec}

\section{Partial traces}

\begin{ssec}
Given a homomorphism 
$\phi :T(A_{1},A_{2},\dots ,A_{l})
 \rightarrow T(A_{1},A_{2},\dots ,A_{l})$,
$b\in \Sigma $ and a basis $\{\epsilon _i\}_{i}$ of $F(b,A_{k})$,
the partial trace of $\phi $ as $A_{k}$ goes
to $b$ is defined to be the morphism
\begin{eqnarray*}
&&{\tt Tr}_{A_{k}\rightarrow b}^{\epsilon _i}(\phi )
=\sum _{i} (id\diamond \epsilon ^{*}_i \diamond id
 )\circ\phi\circ
(id\diamond \epsilon _i \diamond id):\\
&&\qquad\quad T(A_{1},\dots , A_{k-1}, b,A_{k+1}, \dots ,A_{l})\rightarrow
T(A_{1},\dots , A_{k-1}, b,A_{k+1}, \dots ,A_{l}),
\end{eqnarray*}
where 
$\{\epsilon ^{*}_i\}_{i}$ 
is the dual basis of $F(A_{k},b)$, and $(id\diamond \psi\diamond id)$ stands for 
\newline  $T(id_{A_{1}},\dots,id_{A_{k-1}},\psi ,id_{A_{k+1}},\dots 
  ,id_{A_{l}})$. The corresponding diagram is 
presented in figure \ref{diagr.patr} (a). 
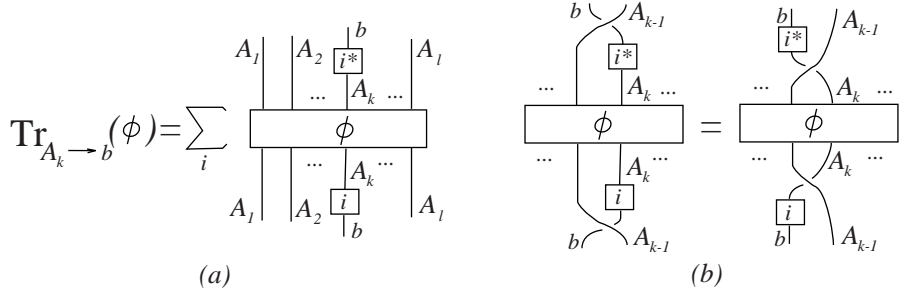
\begin{figure}[h]
\setlength{\unitlength}{1cm}
\begin{center}
\begin{picture}(14,3.5)
\put(0,-0.3){\fig{figptr.pict}}
\end{picture}
\end{center}
\caption{Partial trace.}
\label{diagr.patr}
\end{figure}
\end{ssec}

\begin{ssec} 
\label{ptrprop} 
Let $\phi $ be as above.
\begin{itemize}
\item[(a)] If $\psi:A_{k}\rightarrow A'_{k}$ is an isomorphism and 
  $\epsilon _i,\epsilon _i'$ are bases of $F(b,A_{k})$ and  $F(b,A_{k}')$ 
  correspondingly, then
  $$
  {\tt Tr}_{A_{k}\rightarrow b}^{\epsilon _i}(\phi )={\tt Tr}_{A_{k}'\rightarrow 
  b}^{\epsilon _i'}((id\diamond \psi\diamond id)\circ\phi \circ
  (id\diamond \psi ^{-1}\diamond id));
  $$
\item[(b)] If $k\neq m$, and $a,b\in \Sigma$, then
$${\tt Tr}_{A_{m}\rightarrow a}({\tt Tr}_{A_{k}\rightarrow 
   b}(\phi )) ={\tt Tr}_{A_{k}\rightarrow b}({\tt Tr}_{A_{m}\rightarrow 
   a}(\phi ))\equiv {\tt Tr}_{\{A_{k}\rightarrow b,A_{m}\rightarrow 
   a\}}(\phi );
$$
\item[(c)] If $A_{k}=R(B_{1},B_{2}, \ldots , B_{m})$, 
   i.e. for some $T'$,
   $$
   T(A_{1},A_{2},\dots  ,A_{l})=
   T'(A_{1},A_{2},\dots , A_{k-1}, B_{1}, \ldots ,B_{m},A_{k+1}, \dots 
   ,A_{l}),
$$ 
then 
$$
{\tt Tr}_{A_{k}\rightarrow b}(\phi ) =
 \sum _{a\in \Sigma }{\tt Tr}_
{R(B_{1}, \ldots ,B_{s-1},a, B_{s+1},\ldots , B_{m})\rightarrow b}(
{\tt Tr}_{B_{s}\rightarrow a}(\phi ));
$$
\item[(d)] If $\pi$ is a permutation in ${\bf S}_l$, then 
$
\pi\circ {\tt Tr}_{A_{k}\rightarrow b}(\phi )\circ \pi ^{-1}=
 {\tt Tr}_{A_{\pi (k)}\rightarrow b}
(\pi\circ \phi \circ \pi ^{-1}) .
$
\end{itemize}
\vspace{0.3cm}
(a) follows from the semisimplicity of the category in the following 
way:
\begin{eqnarray*}
&&\sum _i(id\diamond (\epsilon '_i{}^* \circ \psi)\diamond id)\circ \phi\circ 
    (id\diamond (\psi ^{-1} \circ \epsilon '_i )\diamond id)=\\
&&\qquad =\sum _{i,j}(id\diamond (\epsilon '_i{}^* \circ \psi\circ \epsilon _j)
  \diamond id)\circ (id\diamond \epsilon ^*_j \diamond id)\circ \phi\circ 
   (id\diamond (\psi ^{-1} \circ \epsilon '_i )\diamond id)=\\
&&\qquad =\sum _{i,j}(id\diamond \epsilon ^*_j \diamond id)\circ \phi\circ 
   (id\diamond (\psi ^{-1} \circ \epsilon '_i )\diamond id)\circ
    (id\diamond (\epsilon '_i {}^* \circ \psi\circ \epsilon _j)
   \diamond id)=\\
&&\qquad =\sum _j(id\diamond \epsilon ^*_j \diamond id)\circ \phi\circ 
      (id\diamond \epsilon _j \diamond id).
\end{eqnarray*}

(b) follows directly from the definition. To prove (c),
we observe that from \ref{defn.semis} follows that if $\{\tau _{i}\}_{i}$ is
a basis for $F(b,R(B_{1}, \ldots ,B_{s-1},a, B_{s+1},\ldots , B_{m}))$, 
and $\{\epsilon _{j}\}_{j}$ is a basis for $F(a,B_{s})$, then 
$\{R(id_{B_{1}},\ldots ,id_{B_{s-1}},\epsilon _{j},id_{B_{s+1}} \ldots, 
id_{B_{m}})\circ \tau 
_{i}\}_{a, i,j}$ is a basis for $F(b,A_{k})$.  
The use of this basis in the evaluation of the partial trace leads to the 
expression in (c). The statement in (d) is straightforward: every 
permutation is a product of ones of the form $(k-1,k)$, and for those 
is shown in figure \ref{diagr.patr} (b). 

As a consequence of (a) above, the 
indication of the basis in the notation of the partial trace will be 
omitted.
\end{ssec}

\section{Definition of the invariant}

To simplify the notation, we make a standard choice of 
bracketing. 
When $T$ is of the form 
$(\dots ((O\otimes O)\otimes O)\dots \otimes O)\otimes O$,
we will write $T(B^k)$ as simply $B^k$, and $T(B_{1}, B_{2},\dots ,B_{k})$
as $B_{1}\diamond  B_{2}\diamond\dots \diamond B_{k}$. In particular,
$A(BC)B^{3}C=(((A,(BC))((B,B),B))C$.

\begin{ssec}
\label{fbl.prop}
We start by introducing some definitions.
Given any $b\in\Sigma$ define
$$
b^l\equiv \left\{ \begin{array}{cc} b^l\mbox{ if } l> 0, \\
                                   (b^{*})^{|l|}\mbox{ if } l< 0
                 \end{array}\right.
$$
Then if $l>0$, let $f(b,l)=id_{b^{l+1}}:(b,b, \ldots,b)\rightarrow 
(b,b, \ldots,b)$ be the open morphism in \CA , where $In(f(b,l))$ is
the first copy of $b$ in the lower ends and $Out(f(b,l))$ is the last 
copy of $b$ among the upper ends. For $l<0$, we define 
$f(b,l)=f(b^{*},|l|)^{\hat{}}:(b,b^*, \ldots,b^*)\rightarrow (b^*, 
\ldots,b^*,b)\ $. The diagrams corresponding to $f(b,l)$ in both cases are 
presented in figures \ref{fbl} (a) and (b). Note that with this definition,
for any integer $l$, $f(b,l)^{\hat{}}=f(b^*,-l)$. 
\begin{figure}[h]
\setlength{\unitlength}{1cm}
\begin{center}
\begin{picture}(11,3)
\put(0,-0.3){\fig{figfbl.pict}}
\end{picture}
\end{center}
\caption{Diagram presentation of $f(b,l)$.}
\label{fbl}
\end{figure}
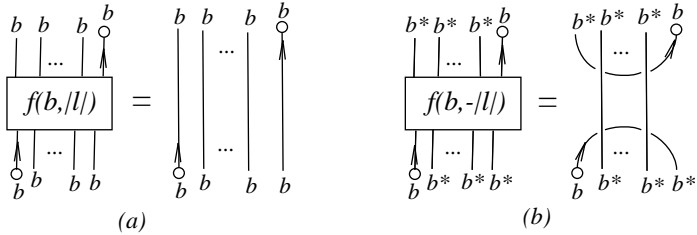
\end{ssec}

\begin{ssec}
\label{fblfbl}
Suppose we are given a nonzero integer $l$ and 
$b_{1},b_{2}, a_{i}\in \Sigma$, $1\leq i\leq |l|$. 
Let
$
\chi :b_{1}b_{1}^lb_{2}b_{2}^l\rightarrow 
(b_{1}b_{2})(\diamond _{i}(b_{1}b_{2})_{i})$  and
$\chi ':b_{1}^lb_{1}b_{2}^lb_{2}\rightarrow 
(b_{1}b_{2})(\diamond _{i}(b_{1}b_{2})_{i})$
be the corresponding permutations, where
$$
(b_{1}b_{2})_{i}=\left\{ \begin{array}{cc} b_{1}b_{2}\mbox{ if } l> 0, \\
                                   b_{1}^{*}b_{2}^{*}\mbox{ if } l< 0
                 \end{array}\right. .
$$
 Then 
${\tt Tr}_{\{(b_{1}b_{2})_{i}\rightarrow a_{i}\}_{i}}
      (\chi '\circ(f(b_{1},l)\otimes f(b_{2},l))\circ\chi ^{-1}):
b_{1}b_{2}(\diamond _{i}a_{i})\rightarrow 
b_{1}b_{2}(\diamond _{i}a_{i})$
is zero unless $a_{i}$ are not all the same, in which case is closely 
related to $f(a_{1},l)$. To be able to give the precise statement,
let 
$\{\epsilon _{k}(a_{i},b_{1}b_{2}):a_{i}\rightarrow b_{1}b_{2}\}_{k}$ be 
a basis for $F(a_{i},b_{1}b_{2})$ and $\{\epsilon 
_{k}(a_{i},b_{1}b_{2}))^*\}_{k}$ be its dual basis. 
 We will denote with the same symbols the open morphism in \CA\ :
$(a_{i})\rightarrow (b_{1},b_{2})$ ($(b_{1},b_{2})\rightarrow (a_{i})$) 
which have a unique 
coupon labeled by $\epsilon _{k}(a_{i},b_{1}b_{2})$ 
(correspondingly $(\epsilon _{k}(a_{i},b_{1}b_{2}))^*$ ), 
and all lower ends are $in$, and all upper ends are $out$. 
Let
$
f^{\epsilon}(a_{i},l)=\sum_{k}\epsilon _{k}(a_{i},b_{1}b_{2})^{*} \circ _{a_{i}}
     f(a_{i},l)\circ _{a_{i}}\epsilon _{k}(a_{i},b_{1}b_{2}):
     b_{1}b_{2}a_{1}^l\rightarrow a_{1}^lb_{1}b_{2}$.
Then if $l>0$ we have 
\begin{eqnarray*}
&&{\tt Tr}_{\{(b_{1}b_{2})_{i}\rightarrow a_{i}\}_{i}}
      (\chi ' \circ(f(b_{1},l)\otimes f(b_{2},l))\circ\chi ^{-1})=
      \sqcap _{i=1}^{l-1} 
      \delta _{a_{i},a_{i+1}}\sum _{k} \pi\circ f^{\epsilon}(a_{1},l).
\end{eqnarray*}
Here $\pi :a_1^lb_1b_2\rightarrow b_1b_2a_1^l$ is the corresponding permutation.
The proof is shown in figure \ref{fptrA}.
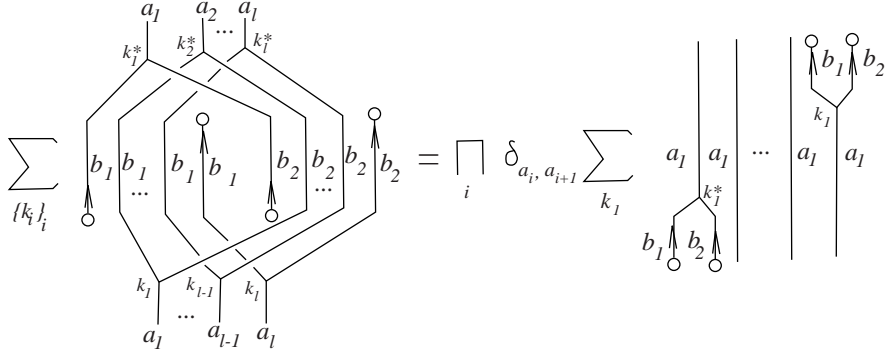
\begin{figure}[h]
\setlength{\unitlength}{1cm}
\begin{center}
\begin{picture}(13,4.5)
\put(0,0){\fig{figpropfbl.pict}}
\end{picture}
\end{center}
\caption{Proof of the property in \ref{fblfbl} for positive $l$.}
\label{fptrA}
\end{figure}

When $l<0$,  the morphisms
$\{
\zeta _{k}(a_{i}^{*},b_{1}^{*}b_{2}^{*})=\gamma _{b_{2}^{*},b_{1}^{*}} \circ 
\hat{\epsilon} _{k} (a_{i},b_{1}b_{2})
\}_{k}$
form a basis for $F(a_{i}^{*},b_{1}^{*}b_{2}^{*})$.
Then taking the partial 
trace with respect to this basis, also in this case, we obtain that
\begin{eqnarray*}
&&{\tt Tr}_{\{(b_{1}b_{2})_{i}\rightarrow a_{i}^{*}\}_{i}}
      (\chi '\circ(f(b_{1},l)\otimes f(b_{2},l))\circ\chi ^{-1})=
      \sqcap _{i=1}^{|l|-1} 
      \delta _{a_{i},a_{i+1}}\sum _{k}\pi\circ f^{\epsilon}(a_{1},l).
\end{eqnarray*}
The proof is shown in figure \ref{fptrB}.
\begin{figure}[h]
\setlength{\unitlength}{1cm}
\begin{center}
\begin{picture}(13,4.5)
\put(0,0){\fig{figprodfbl1.pict}}
\end{picture}
\end{center}
\caption{Proof of the property in \ref{fblfbl} for negative $l$.}
\label{fptrB}
\end{figure}
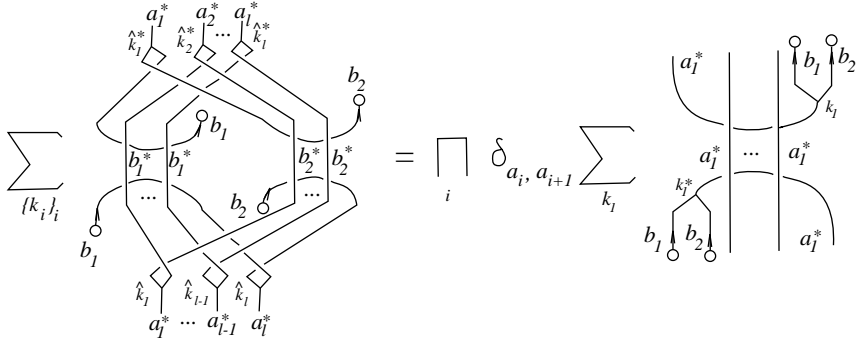
\end{ssec}

\begin{ssec}
First we will define the invariant for the case of a presentation 
 with one relation $P=<x_1,x_2,\dots ,x_n|\; R>$, where
$$
R=x_{i_{1}}^{l_{1}}x_{i_{2}}^{l_{2}} \ldots x_{i_{k}}^{l_{k}}.
$$
To the relation $R$ and to any $b\in\Sigma$ we associate an element
 $R(b)=b^{l_{1}}b^{l_{2}} \ldots b^{l_{k}}\in\A $ and 
a morphism $[R,b]:R(b)\rightarrow R(b)$ defined in the following 
way: 
$$
[R,b]=\circ _{b}(r(b)\,f(b,l_{1})\circ _{b}f(b,l_{2})\circ _{b}\ldots
      \circ _{b} f(b,l_{k})).
$$
In other words, the $out$-end of $f(b,l_{i})$ is connected with the 
$in$-end of $f(b,l_{i+1})$ and at the end the morphism 
is closed by connecting the $out$-end of
$f(b,l_{k})$ with the $in$-end of $f(b,l_{1})$. The result is multiplied
by the rank of $b$. The open morphism 
without the final closure will be denoted with $[R,b]^{o}$.
We will refer to the factors $b^l$ in $R(l)$ as
$b$-factors. 
As  examples, the diagrams corresponding to $[x^2y^{-2}x^{-1}y,b]^o$
and $[x^2y^{-2}x^{-1}y,b]$ are presented in figure \ref{exampl.diag}
(a) and (b).
\end{ssec}
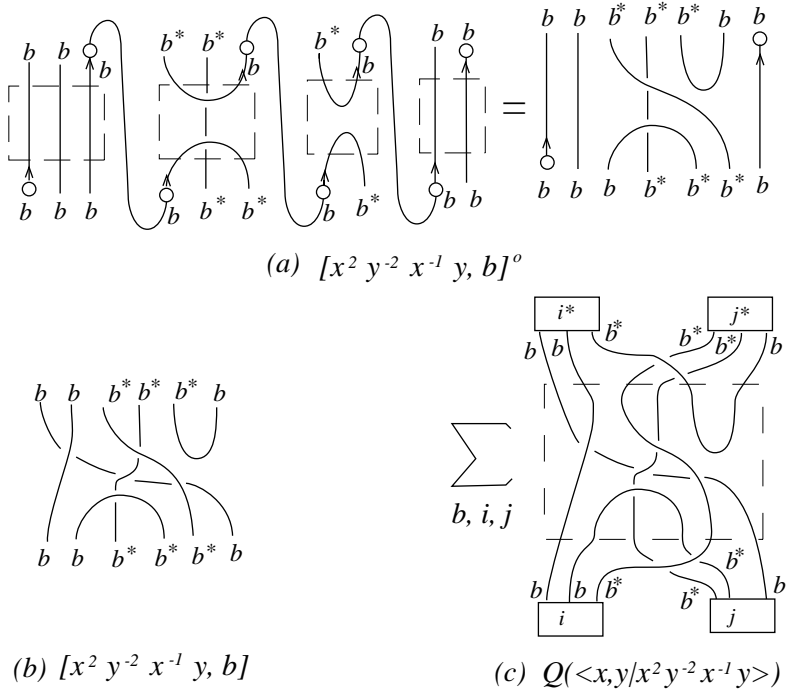
\begin{figure}[h]
\setlength{\unitlength}{1cm}
\begin{center}
\begin{picture}(15,9)
\put(0,-0.3){\fig{figgeninv.pict}}
\end{picture}
\end{center}
\caption{Diagram presentation of $Q(<x,y|x^2y^{-2}x^{-1}yx>)$.}
\label{exampl.diag}
\end{figure}

\begin{ssec}
\label{onerel}
Now to the relation $R$, to $b\in\Sigma $, and to 
any generator $x_{k}$, we associate
an element $g_{k}(R,b)\in\A$, defined as the product of all b-factors 
in $R(b)$ corresponding to the generator $x_{k}:
g_{k}(R,b)=\diamond _{r|x_{i_{r}}=x_{k}}b^{l_{r}}$.
Let $\kappa (R):R(b)\rightarrow g_{1}(R,b)g_{2}(R,b) \ldots 
g_{n}(R,b)$ be the permutation of factors. Then 
the invariant of the presentation $P$ is defined to be:
$$
Q(P)=\sum _{b\in\Sigma }{\tt 
Tr}_{\{g_{k}(R,b)\rightarrow \1\}_{k}}(\kappa (R)\circ[R, b]
       \circ \kappa (R)^{-1}).
$$
The diagram corresponding to 
$Q(<x,y|x^2y^{-2}x^{-1}y>)$ is presented in figure 
\ref{exampl.diag} (c).
The evaluation of the invariant for some simple examples is presented on
figures \ref{examples} and the results are collected below.
\begin{itemize}
\item[(a)] $Q(<x|x^k>)=\sum _{b\in\Sigma}r(b){\tt Trace}
      [\rho _{\1}[b^k](1,2,3,\ldots ,k)]$. In particular,\newline 
$
Q(<x|x^2>)=\sum _{b\in\Sigma |b=b^{*}}r(b);
$
\item[(b)] The invariant for the standard presentation of the fundamental
    group of a surface of genus $n$ is $Q(<x_1,x_2,\dots ,x_{2n} |\;\sqcap 
    _{j=1}^{n-1}x_{2j}x_{2j+1}x_{2j}^{-1}x_{2j+1}^{-1}>)= \sum _{b\in\Sigma 
                                          }r(b)^{-2(n-1)}$;
\item[(c)] $Q(<x,y| xyxy^{-1}x^{-1}y^{-1}>)=\sum _{b\in\Sigma 
      |b=b^{*}}{\tt Trace}[\rho _{\1}[b^3](1,2,3)]$;
\item[(d)] $Q(<x,y| xyx^{-1}y>)=|\{b\in\Sigma 
                                  |b=b^{*}\}|$.
\end{itemize}
Here ${\tt Trace}$ denotes the usual trace of the corresponding 
element in the representation of the symmetric group.
\end{ssec}
\begin{figure}[h]
\setlength{\unitlength}{1cm}
\begin{center}
\begin{picture}(15,13)
\put(0,-0.3){\fig{figQs.pict}}
\end{picture}
\end{center}
\caption{Some examples.}
\label{examples}
\end{figure}

\begin{ssec}
\label{morerel}
Now we extend the definition to the case of presentations with more then one
relation. Let $P=<x_1,x_2,\dots ,x_n|\; R_1,R_2,\dots ,R_m>$ be such a
presentation,
and let $\underline{b}=(b_{1},b_{2} \ldots b_{m})\in \Sigma ^{\times m}$. 
Then  define
\begin{eqnarray*}
&&Rel(\underline{b})=\diamond _{j=1}^m R(b_{j}),Ê\\
&&[P,\underline{b}]=[R_{1},b_{1}]\diamond [R_{2},b_{2}]\diamond \ldots
 \diamond [R_{m},b_{m}]:Rel(\underline{b})\rightarrow Rel(\underline{b}),\\
&&G_{k}(\underline{b})=\diamond _{j=1}^m 
        g_{k}(R_{1},b_{1})g_{k}(R_{2},b_{2}) \ldots 
                       g_{k}(R_{m},b_{m}),\\
&&\xi (P): Rel(\underline{b})\rightarrow \diamond _{k=1}^nG_{k}(\underline{b})
    \mbox{ to be the corresponding permutation of $b$-factors}.
\end{eqnarray*}
Note that the permutation $\xi (P)$ groups the 
$b$-factors in  $Rel(\underline{b})$ according to the generator to which 
they correspond. Then  
the invariant of the presentation $P$ is defined to be
$$
Q(P)=\sum _{\underline{b}\in\Sigma ^{\times m}}{\tt 
Tr}_{\{G_{k}(\underline{b})\rightarrow \1\}_{k}}
        (\xi(P)\circ[P, \underline{b}]
       \circ \xi (P)^{-1}).
$$
As example the diagram corresponding to $Q(<x,y|x^2y^2, 
xyx^{-1}y^{-1}>)$ is shown in figure \ref{tworel}.
\begin{figure}[h]
\setlength{\unitlength}{1cm}
\begin{center}
\begin{picture}(12,4)
\put(0,0){\fig{figtworel.pict}}
\end{picture}
\end{center}
\caption{Example with two relations.}
\label{tworel}
\end{figure}

\begin{theo}
$Q(P)$ is invariant under the AC-moves (i)-(vi).
\end{theo}
The proof is contained in the next section.
\end{ssec}

\section{Invariance under the AC - moves}

We consider each move separately. The presentation obtained after 
performing the move will be denoted with $P'$ and the corresponding
products of $b$-factors will be denoted 
with $Rel'(\underline{b})$ and $G_{k}'(\underline{b})$.

\begin{ssec}
The first move consists in exchanging the places of $R_{1}$ and 
$R_{s}$. Let 
$$\underline{b}'=(b_{s}, b_{2},b_{3}, \ldots , 
b_{s-1},b_{1},b_{s+1}, \ldots, b_{m})$$ and let 
$\pi _{k}:G_{k}(\underline{b})\rightarrow 
G_{k}'(\underline{b}')$ be the permutation which exchanges the places 
of $g_{k}(R_{1},b_{1})$ and $g_{k}(R_{s},b_{s})$. Then  $\pi = \pi 
_{1}\diamond \pi _{2}\diamond \ldots \diamond \pi _{n}: 
\diamond _{k=1}^nG_{k}(\underline{b})\rightarrow 
\diamond _{k=1}^nG_{k}'(\underline{b}')$. As morphisms in \CA\ we 
have that
$$
\pi\circ \xi (P)\circ [P,\underline{b}]\circ \xi (P)^{-1}\circ \pi 
^{-1}= \xi (P')\circ [P',\underline{b}']\circ \xi (P')^{-1}.
$$
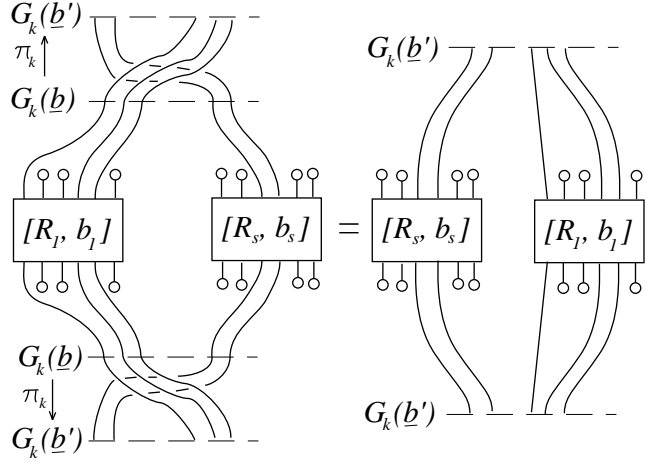
\begin{figure}[h]
\setlength{\unitlength}{1cm}
\begin{center}
\begin{picture}(9,6)
\put(0,0){\fig{figmove1.pict}}
\end{picture}
\end{center}
\caption{Invariance under move (i).}
\label{moveone}
\end{figure}
This fact is illustrated in figure \ref{moveone}. Then according to
\ref{ptrprop} (a) we have
\begin{eqnarray*}
&Q(P')&=\sum _{\underline{b}'\in\Sigma ^{\times m}}{\tt 
    Tr}_{\{G_{k}'(\underline{b}')\rightarrow \1\}_{k}}
        (\pi\circ\xi(P)\circ [P, \underline{b}]
       \circ \xi (P)^{-1}\circ\pi ^{-1}) = \\
&&=\sum _{\underline{b}'\in\Sigma ^{\times m}}{\tt 
    Tr}_{\{G_{k}(\underline{b})\rightarrow \1\}_{k}}
        (\xi(P)\circ[P, \underline{b}]
       \circ \xi (P)^{-1}) = Q(P).
\end{eqnarray*}
This shows the invariance under the first move.
\end{ssec}

\begin{ssec}
The second move conjugates $R_{1}$ with an arbitrary element in the 
group. Obviously, it is enough to show the invariance under conjugation
with one of the generators $R_{1}\rightarrow R_{1} '=x_{l}R_{1}x_{l}^{-1}$. 
Then $G_{k}$ changes only when $k=l$, and 
$G_{l}'(\underline{b})= b_{1}g_{l}(R_{1},b_{1})b_{1}^*g_{l}(R_{2},b_{2})\ldots
g_{l}(R_{m},b_{m})$. Let 
$$
\pi :(\diamond _{k<l}G_{k}(\underline{b}))
G_{l}'(\underline{b})
(\diamond _{k>l}G_{k}(\underline{b}))
\rightarrow (\diamond _{k<l}G_{k}(\underline{b})) 
b_{1}b_{1}^* G_{l}(\underline{b})
(\diamond _{k>l}G_{k}(\underline{b}))
$$
be the corresponding permutation of factors. The goal now is to show 
that 
$$
E={\tt Tr}_{G_{l}'(\underline{b})\rightarrow \1}
        (\xi(P')\circ[P', \underline{b}]\circ\xi (P')^{-1})=
{\tt Tr}_{G_{l}(\underline{b})\rightarrow \1}
        (\xi(P)\circ[P, \underline{b}]\circ\xi (P)^{-1}),
$$ 
which according to \ref{ptrprop} (b) would imply the statement. 
By using properties \ref{ptrprop} (a) and (c) of the partial trace 
we obtain
\begin{eqnarray*}
&E&={\tt Tr}_{b_{1}b_{1}^*G_{l}(\underline{b})\rightarrow \1}
   (\pi \circ\xi(P')\circ[P', \underline{b}]\circ\xi (P')^{-1}\circ\pi ^{-1})=\\
&&=\sum _{a\in\Sigma}{\tt Tr}_{aG_{l}(\underline{b})\rightarrow \1}[
     {\tt Tr}_{b_{1}b_{1}^*\rightarrow a}
 (\pi \circ\xi(P')\circ[P', \underline{b}]\circ\xi (P')^{-1}\circ\pi ^{-1})].
\end{eqnarray*}
We remind that $[P',\underline{b}]=\diamond _{i=1}^{m}[R'_{i},b_{i}]$, 
but
the partial trace  ${\tt Tr}_{b_{1}b_{1}^*\rightarrow a}$ actually 
involves only the morphism $[R'_{1},b_{1}]$, in sense that 
$[R'_{1},b_{1}]$ maps $b_1R_1(b_1)b_1^*$ into itself, as shown on
figure \ref{movetwo} (a), and the effect of
the permutation $\pi$ is to move the $b$-factors $b_1$ and $b_1^*$ together,
so that the partial trace can be taken. Then 
figure \ref{movetwo} (b)  proves that
$$
{\tt Tr}_{b_{1}b_{1}^*\rightarrow a}
 (\pi \circ\xi(P')\circ[P', \underline{b}]\circ\xi (P')^{-1}\circ\pi ^{-1})=
      \delta _{a,{\bf 1}}\,\nu \circ(id_{{\bf 1}}\diamond (\xi(P)
      \circ [P, \underline{b}]
      \circ\xi (P)^{-1}))\circ\nu ^{-1}.
$$
Here $\nu: \1\, G_{1}(\underline{b})G_{2}(\underline{b}) \ldots 
G_{2}(\underline{b})\rightarrow G_{1}(\underline{b})G_{2}(\underline{b}) \ldots 
G_{l-1}(\underline{b}) \,\1\, G_{l}(\underline{b}) \ldots
G_{n}(\underline{b})$ is the corresponding permutation. The factor 
\1\ can be removed from ${\tt Tr}_{\1 G_{l}(\underline{b})\rightarrow 
\1}$ by conjugating the argument with $\etaa{G_{l}(\underline{b})}$ 
and using \ref{ptrprop} (a). 
This completes the proof of the invariance under (ii).
\begin{figure}[h]
\setlength{\unitlength}{1cm}
\begin{center}
\begin{picture}(10,7.3)
\put(0,0){\fig{figmove2.pict}}
\end{picture}
\end{center}
\caption{Invariance under move (ii).}
\label{movetwo}
\end{figure}
\end{ssec}

\begin{ssec}
The third move consists of  replacing $R_{1}$ with $R_{1}^{-1}$. First
we observe that if $R=x_{i_{1}}^{l_{1}}x_{i_{2}}^{l_{2}} \ldots 
x_{i_{s}}^{l_{s}}$ is a relation, then
\begin{eqnarray*}
&[R^{-1},b]&=r(b)\circ _{b}(f(b,-l_{s})\circ _{b}f(b,-l_{s-1})\circ _{b}
\ldots \circ _{b}f(b,-l_{1}))=\\
&&=r(b^{*})\circ _{b}(f(b^{*},l_{s})^{\hat{}}\circ _{b}
    f(b^*,l_{s-1})^{\hat{}}\circ _{b}
\ldots \circ _{b}f(b^*,l_{1})^{\hat{}}\,)=\\
&&=\pi (R)\circ [R,b^*]\circ \pi (R)^{-1},
\end{eqnarray*}
where $\pi (R):b^{-l_{1}}b^{-l_{2}} \ldots b^{-l_{s}}\rightarrow 
          b^{-l_{s}}b^{-l_{s-1}} \ldots b^{-l_{1}}$ is the permutation of
$b$-factors. The last equality is illustrated in figure \ref{movethree}.
\begin{figure}[h]
\setlength{\unitlength}{1cm}
\begin{center}
\begin{picture}(12,4)
\put(0,-0.3){\fig{figmove3.pict}}
\end{picture}
\end{center}
\caption{Invariance under move (iii).}
\label{movethree}
\end{figure}
In the case we want to study, $\pi (R_{1})$ induces a permutation 
$\pi = \pi (R_{1})\diamond id :R_{1}(b_{1}^{*})\diamond (\diamond 
_{i=2}^{m}R_{i}(b_{i}))\rightarrow R_{1}^{-1}(b_{1})\diamond (\diamond 
_{i=2}^{m}R_{i}(b_{i}))$. 
Then from the discussion above it follows that 
$$
[P',\underline{b}]=\pi \circ [P,\hat{\underline{b}}]\circ \pi ^{-1},
$$
where $\hat{\underline{b}}=(b_{1}^{*}, b_{2}, \ldots, b_{m})$.
Moreover, $G_{k}'(\underline{b})$ can be obtained from 
$G_{k}(\hat{\underline{b}})$ 
by reversing the order of the $b$-factors of the form $b^{l_{j}}_{1}$. 
Let the corresponding permutation be 
$\chi _{k}: G_{k}'(\underline{b})\rightarrow 
G_{k}(\hat{\underline{b}})$, and let $\chi =\chi _{1}\diamond \chi 
_{2}\diamond \ldots \diamond \chi _{m}: \diamond _{k=1}^n
G_{k}'(\underline{b})\rightarrow \diamond _{k=1}^n G_{k}(\hat{\underline{b}})$.
Then from property \ref{ptrprop} we obtain
\begin{eqnarray*}
&Q(P')&=\sum _{\underline{b}\in\Sigma ^{\times m}}{\tt 
          Tr}_{\{G_{k}'(\underline{b})\rightarrow \1\}_{k}}
        (\xi(P')\circ[P', \underline{b}]
       \circ \xi (P')^{-1})=\\
&&=\sum _{\underline{b}\in\Sigma ^{\times m}}{\tt 
          Tr}_{\{G_{k}(\hat{\underline{b}})\rightarrow \1\}_{k}}
        (\chi\circ \xi(P')\circ
        \pi \circ [P,\hat{\underline{b}}]\circ \pi ^{-1}
       \circ \xi (P')^{-1}\circ \chi ^{-1})=\\
&&=\sum _{\underline{b}\in\Sigma ^{\times m}}{\tt 
          Tr}_{\{G_{k}(\hat{\underline{b}})\rightarrow \1\}_{k}}
        (\xi(P)\circ[P, \hat{\underline{b}}]
       \circ \xi (P)^{-1})=Q(P).
\end{eqnarray*}
This completes the proof of the invariance under move (iii).
\end{ssec}

\begin{ssec}
The forth move consists of replacing $R_{1}$ with $R_{1}R_{2}$. Figure
\ref{movefourA} presents the corresponding change in $[P,\underline{b}]$.
\begin{figure}[h]
\setlength{\unitlength}{1cm}
\begin{center}
\begin{picture}(12,2)
\put(0,0){\fig{figmove4.pict}}
\end{picture}
\end{center}
\caption{Change in $[P,\underline{b}]$ under move (iv).}
\label{movefourA}
\end{figure} 
Let 
$
R_{2}=x_{i_{1}}^{l_{1}}x_{i_{2}}^{l_{2}} \ldots 
x_{i_{s}}^{l_{s}}
$,
and \newline
$\chi : R_{1}(b_{1})R_{2}(b_{1})R_{2}(b_{2})\rightarrow
R_{1}(b_{1})\diamond _{i=1}^s \diamond _{j=1}^{|l_{i}|}(b_{1}b_{2})_{i,j}$ 
be the corresponding permutation, where
$$
(b_{1}b_{2})_{i,j}=\left\{\begin{array}{c} b_{1}b_{2} \mbox{ if } 
                                                    l_{i}>0\\
                                   b_{1}^{*}b_{2}^{*} \mbox{ if } 
                                                    l_{i}<0
                        \end{array}\right. .
$$
Then $\chi \circ([R_{1}R_{2},b_{1}]\diamond [R_{2},b_{2}])\circ \chi ^{-1}$ maps  
$R_{1}(b_{1})\diamond _{i=1}^s \diamond _{j=1}^{|l_{i}|}
  (b_{1}b_{2})_{i,j}$ into itself. As a first step we prove that
\begin{eqnarray*}
&W&=\sum _{b_{2}\in\Sigma} {\tt Tr}_{\{(b_{1}b_{2})_{i,j}
\rightarrow a_{i,j}\}_{i,j}}
    (\chi \circ([R_{1}R_{2},b_{1}]\diamond [R_{2},b_{2}])\circ\chi ^{-1})=\\
&&\qquad \qquad   \sqcap _{i=1}^{s-1}(
      \delta _{a_{i,1},a_{i+1,1}}
      \sqcap _{j=1}^{|l_{i}|-1} 
      \delta _{a_{i,j},a_{i,j+1}})
    ([R_{1},b_{1}]\diamond [R_{2},a_{1,1}])
\end{eqnarray*}
Let $a_{i}\equiv a_{i,1}$, $1\leq i\leq s$. Recalling
 \ref{fbl.prop}, where the partial trace involving $f(b_1,l)\otimes f(b_2,l)$
was evaluated, and we have that $W=\circ 
_{b_{1}}([R_{1},b_{1}]^{o}\circ _{b_{1}}G)$, where
$$
G=(\sqcap _{i=1}^{s-1}\sqcap _{j=1}^{|l_{i}|-1} 
      \delta _{a_{i,j},a_{i,j+1}})\sum _{b_{2}\in\Sigma} \circ _{b_{2}}
 (f^{\epsilon}(a_{1},l_{1}) \circ _{\{b_{1},b_{2}\}}f^{\epsilon}(a_{2},l_{2})
 \circ _{\{b_{1},b_{2}\}}
    \ldots \circ _{\{b_{1},b_{2}\}}f^{\epsilon}(a_{s},l_{s})).
$$
Then the statement for $W$ follows from the fact that 
$\epsilon _{k}(a_{i},b_{1}b_{2})^{*}\circ 
\epsilon _{l}(a_{i+1},b_{1}b_{2})=\delta _{a_{i},a_{i+1}}\delta 
_{k,l}\,id_{a_{i}}$ and from proposition \ref{dbasis} as it is shown in figure
\ref{movefourB}. 
\begin{figure}[h]
\setlength{\unitlength}{1cm}
\begin{center}
\begin{picture}(9,6)
\put(0,-0.3){\fig{figmove41.pict}}
\end{picture}
\end{center}
\caption{Invariance under move (iv).}
\label{movefourB}
\end{figure} 

Now the invariance under move (iv) is easy to show. Fix $k$ between
1 and $n$.
Then after the move, 
$G_{k}'(\underline{b})=g_{k}(R_{1},b_{1})g_{k}(R_{2},b_{1})
g_{k}(R_{2},b_{2})\ldots g_{k}(R_{m},b_{m})$. Let
$G^{>2}_{k}(\underline{b})=\diamond _{i=3}^s g_{k}(R_{i},b_{i})$ and
 $g_{k}(R_{2},b)=\diamond _{r=1}^{s_{k}}b^{l_{i_{r}}}$.
Then define $\tau _{k}: G_{k}'(\underline{b})\rightarrow 
g_{k}(R_{1},b_{1})\diamond _{r=1}^{s_{k}} 
\diamond _{j=1}^{|l_{i_{r}}|}
(b_{1}b_{2})_{i_{r},j}G^{>2}_{k}(\underline{b})$ to be the corresponding 
permutation, and 
let $\tau =\tau _{1}\diamond \tau _{2}\diamond \ldots \diamond \tau _{m}:
\diamond _{k=1}^{n}G_{k}'(\underline{b})\rightarrow 
\diamond _{k=1}^{n} g_{k}(R_{1},b_{1})\diamond _{r=1}^{s_{k}}
\diamond _{j=1}^{|l_{i_{r}}|}
(b_{1}b_{2})_{i_{r},j}G^{>2}_{k}(\underline{b})$. Then
from \ref{ptrprop} (a) and (c) we obtain,
\begin{eqnarray*}
&&\sum _{\underline{b}\in\Sigma ^{\times m}}{\tt 
          Tr}_{G_{k}'(\underline{b})\rightarrow \1 }
        (\xi(P')\circ[P', \underline{b}]
       \circ \xi (P')^{-1})=\\
&&\quad =\sum _{\underline{b}\in\Sigma ^{\times m}}{\tt 
          Tr}_{g_{k}(R_{1},b_{1})\diamond _{r=1}^{s_{k}} 
        \diamond _{j=1}^{|l_{i_{r}}|}
         (b_{1}b_{2})_{i_{r},j}G^{>2}_{k}(\underline{b})\rightarrow \1}
        (\tau \circ \xi (P')\circ[P', \underline{b}]
       \circ \xi (P')^{-1}\circ \tau ^{-1})=\\
&&\quad =\sum _{\underline{b}\in\Sigma ^{\times m}}
   \sum _{\underline{a}_{r}\in\Sigma ^{\times s_{k}}}{\tt 
          Tr}_{g_{k}(R_{1},b_{1})\diamond _{r=1}^{s_{k}} 
        \diamond _{j=1}^{|l_{i_{r}}|}
          a_{r,j}G^{>2}_{k}(\underline{b})\rightarrow \1}{\tt 
          Tr}_{\{(b_{1}b_{2})_{i_{r},j}\rightarrow a_{r,j}\}_{r,j}}\\
&&\qquad\qquad\qquad   (\tau \circ \xi(P')\circ[P', \underline{b}]
       \circ \xi (P')^{-1}\circ \tau ^{-1}).
\end{eqnarray*}
Now we observe that
$$
\tau \circ \xi(P')\circ[P', \underline{b}]
       \circ \xi (P')^{-1}\circ \tau ^{-1}=
\xi (P)\circ (
(\chi \circ([R_{1}R_{2},b_{1}]\diamond [R_{2},b_{2}])\circ \chi ^{-1})
\diamond _{i=3}^m [R_i,b_i])\circ \xi (P) ^{-1}.
$$
Here the use of $\xi (P)$ is somewhat abusive. By definition
$\xi (P)$  is the permutation of $b$-factors 
$\diamond _{j=1}^m R(b_j)\rightarrow
\diamond _{k=1}^{n}G_k(\underline{b})$, 
while on the r.h.s. above
the same notation is used to indicate the same
permutation between expressions where each factors $b_2^{l_i}$ has 
been replaced with $\diamond _{j=1}^{l_i}(b_1b_2)_{i,j}$. 
Then the statement follows from \ref{ptrprop} (d).
\end{ssec}

\begin{ssec}
The invariance under the moves (v) and (vi) is straightforward: according 
to the definition, these moves change the value of the invariant by
multiplication or division by
$$
\sum _{b\in\Sigma}r(b){\tt Tr}_{b\rightarrow \1}id _b =1.
$$
\end{ssec}

\section{The definition of Quinn's invariant}

\begin{ssec}
In this section we show that the invariant defined above is actually the
one produced by the algorithm described in  \cite{Q:lectures}. 
First we introduce some morphisms  which are being
used in the algorithm.
\begin{itemize}
\item[(i)] Given any $a,b,c\in\Sigma $ define
$cycl(a,b,c) :F(a,b c)\rightarrow F(b^*,c a^*)$ as
\begin{eqnarray*}
&cycl(a,b,c)(\phi ):&b^{*}\rightarrow b^{*}\1 \mapright{\Lambda _{a^{*}}}
b^{*}(aa^{*})\mapright{id_{b^{*}}\diamond (\phi\diamond id_{a^*})}\\
&&b^{*}((bc)a^{*})\rightarrow b^{*}b(ca^{*})
\mapright{\lambda _{b^{*}}}\1 (ca^{*})\rightarrow ca^{*}.
\end{eqnarray*}
$cycl(a,b,c)$ is actually an isomorphisms with inverse $
cycl(c^*,a^*,b)\circ cycl(b^*,c,a^*)$.
\item[(ii)] Given any $a,b\in\Sigma $ and $A\in\Sigma$ 
define
\begin{eqnarray*}
&m_{b}:F(a,A)\longrightarrow F(\1,(b b^*) A),
    & \psi\rightarrow (\Lambda _{b^*}\diamond \psi )\circ\etaa{a}^{-1} \\
&d_{b}:F(\1,(b b^*) A) \longrightarrow F(a,A),
   &\phi \rightarrow \etaa{A}\circ (\lambda _b\diamond id_A)\circ 
   \phi.
\end{eqnarray*}
\end{itemize}
The corresponding  diagrams are presented in figures 
\ref{diagr.md}.
\begin{figure}[h]
\setlength{\unitlength}{1cm}
\begin{center}
\begin{picture}(12,3)
\put(0,0){\fig{figmd.pict}}
\end{picture}
\end{center}
\caption{The morphisms $m_{b}$ and $d_{b}$.}
\label{diagr.md}
\end{figure}
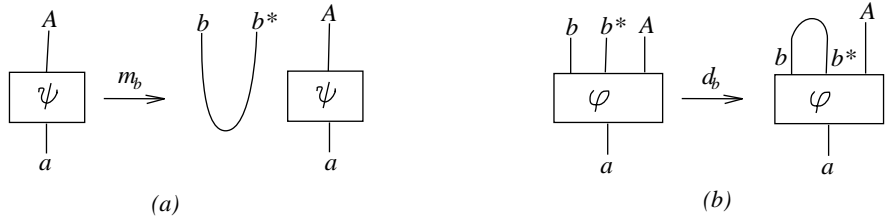 
\end{ssec}

\begin{ssec}
\label{algorithm}
The invariant of a presentation
$P=(x_1,x_2,\dots ,x_n|\, R_1,R_2,\dots, R_m)$
as defined in \cite{Q:lectures} is  a map
$Q(P):F(\1,\1)\rightarrow F(\1,\1)$, i.e. an element in $F(\1,\1)\simeq 
K$. This map is obtained as a composition of the
morphisms which are listed below. Here we will refer to the space 
$V(k)=\oplus _{y_j\in\Sigma} F(\1,y_1y^*_1y_2y^*_2\dots y_ky^*_k)$
as the {\em state space 
of $k$ generators}, 
and to the summand $F(\1,\1\,\1\dots \1\,\1)$ of $V(k)$ as the 
trivial summand of this state space.

\begin{itemize}
\item[(i)] {\bf Beginning presentation.} This morphism maps $F(\1,\1)$ into
$V(n)$ by embedding $F(\1,\1)$ into the trivial summand.
\begin{eqnarray*}
K=&F(\1,\1)\rightarrow   &\oplus _{y_j\in\Sigma} 
               F(\1,y_1y^*_1 y_2y^*_2\dots y_ny^*_n),\\
    &1     \rightarrow   &1\in F(\1,\1\,\1\dots \1\,\1).
\end{eqnarray*}

\item[(ii)] {\bf Beginning relation.} Starting a relation leads to
the appearance of an additional generator, and therefore maps $V(n)$
into $V(n+1)$ by using the map $m_b$:
$$
\oplus _{y_j\in\Sigma} F(\1,y_1y^*_1y_2y^*_2\dots y_ny^*_n) 
\stackrel{\oplus _br(b)^2m_{b^*}}{\longrightarrow}
\oplus _{b,y_j\in\Sigma}
     F(\1,b^*by_1y^*_1y_2y^*_2\dots y_ny^*_n).
$$

\item[(iii)] {\bf The morphism
corresponding to a factor $(y_i)^s$ in a relation.} 
The main ingredient in it is the {\em circulator} 
which is the following map of \newline 
$\oplus _{x\in\Sigma }F(a,bxx^*)$ in itself.
\begin{eqnarray*}
&CR(a,b): &\oplus _{x\in\Sigma }F(a,bxx^*)
        \stackrel{(1,3,2)}{\longrightarrow} 
          \oplus _{x\in\Sigma }F(a,x(x^*b))
      \stackrel{\nabla ^{-1}}{\rightarrow}\\
    &&\oplus _{x,z\in\Sigma } F(a,xz^*)\otimes F(z^*,x^*b)
      \stackrel{id\diamond cycl}{\longrightarrow} 
      \oplus _{x,z\in\Sigma } F(a,xz^*)\otimes F(x,bz) 
      \stackrel{\nabla }{\rightarrow}\\
    && \oplus _{z\in\Sigma }F(a,bzz^*). 
\end{eqnarray*}
Here $\nabla $ is as in \ref{defn.semis}. The diagram describing $CR(a,b)$ is
presented in figure \ref{def.circ}. We note that the circulator is actually
an isomorphism and a specific expression for its inverse will
be provided later. 
\begin{figure}[h]
\setlength{\unitlength}{1cm}
\begin{center}
\begin{picture}(5.5,3.7)
\put(0,-0.3){\fig{figbascirc.pict}}
\end{picture}
\end{center}
\caption{The circulator.}
\label{def.circ}
\end{figure}
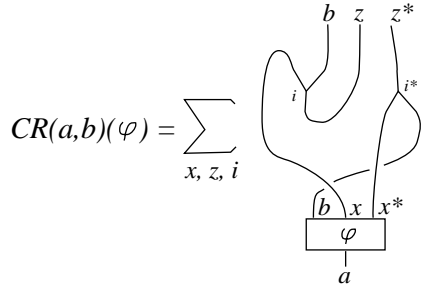 
Then a factor $(y_i)^s$ in a relation corresponds to the
following composition of morphisms. First, in the
state space of the $n+1$ generators,  $y_iy^*_i$ is moved
next to $b$. Then the part corresponding to
$by_iy^*_i$ is separated from the rest by using the isomorphism $\nabla ^{-1}$,
and the circulator $CR$ is applied $s$ times on it. Then following
the inverse steps one goes back to the original form
of the state space. 
\begin{eqnarray*}
&&\oplus _{b,y_j\in\Sigma} F(\1,b^*by_1y^*_1
       y_2y^*_2\dots y_ny^*_n)
   \rightarrow \\
&&\oplus _{b,y_j\in\Sigma} F(\1,b^*(by_iy^*_i)y_1y^*_1\dots
   \hat{y}_i\hat{y^*}_i\dots y_ny^*_n) 
   \stackrel{\nabla ^{-1}}{\rightarrow } \\
&&\oplus _{y_j\in\Sigma ,j\neq i}\oplus _{a,b\in\Sigma} F(\1,b^*ay_1y^*_1\dots
  \hat{y}_i\hat{y^*}_i\dots y_ny^*_n) \otimes 
  (\oplus _{y_i\in\Sigma }F(a,by_iy^*_i))
   \stackrel{id\otimes (CR(b,a))^s}{\longrightarrow}\\
&&\oplus _{y_j\in\Sigma ,j\neq i}\oplus _{a,b\in\Sigma} F(\1,b^*ay_1y^*_1\dots
   \hat{y}_i\hat{y^*}_i\dots y_ny^*_n) \otimes 
   (\oplus _{y_i\in\Sigma }F(a,by_iy^*_i))
   \stackrel{\nabla}{\rightarrow } \\
&&\oplus _{b,y_j\in\Sigma} F(\1,b^*(by_iy^*_i)y_1y^*_1\dots
  \hat{y}_i\hat{y^*}_i\dots y_ny^*_n) 
  \rightarrow 
\oplus _{b,y_j\in\Sigma} F(\1,b^*by_1y^*_1
  y_2y^*_2\dots y_ny^*_n),
\end{eqnarray*}
where $\hat{y}$ indicates that the corresponding term is missing.

\item[(iv)] {\bf Ending relation.} This is a map from 
$V(n+1)$ generators back into $V(n)$  proportional
to the one induced by $d_b$.
$$
\oplus _{b,y_j\in\Sigma}
     F(\1,b^*by_1y^*_1y_2y^*_2\dots y_ny^*_n)
\stackrel{\oplus _br(b)^{-1}d_b}{\longrightarrow}
\oplus _{y_j\in\Sigma} F(\1,y_1y^*_1y_2y^*_2\dots y_ny^*_n). 
$$

\item[(v)] {\bf Ending presentation.}
              This is a map from $V(n)$ generators back
     to the ring $V(0)$, which is injective on the trivial summand on the
     state space and sends all other summands of $V(n)$ to 0.
$$
\oplus _{y_j\in\Sigma} 
               F(\1,y_1y^*_1 y_2y^*_2\dots y_ny^*_n)
\rightarrow F(\1,\1)   ,\;
    1\in F(\1,\1\,\1\dots \1\,\1)\rightarrow 1.        
$$
\end{itemize}
\end{ssec}

\begin{theo} \label{qequiv} Given a group presentation $P$, the invariant 
      of $P$, produced by the above algorithm, is exactly $Q(P)$. 
\end{theo}
The key is to understand the powers of the circulator. Given 
$b\in \Sigma$ and a nonzero integer $l$,  we use the definitions of $b^l$
and $f(b,l)$ made in \ref{fbl.prop}. Let $x,w\in \Sigma$, 
$\{\epsilon _i\}_{i}$  be a basis for
$F(w^*,x^*b^l)$ and let $\{\epsilon _i^{*}\}_{i}$ be its dual. 
We define morphisms 
$\phi (l,b,x,w,i): bxx^*\rightarrow x^*(bb^l)w$ and 
$\psi (l,b,x,w,i): x^*(b^lb)w\rightarrow bww^*$ in the following way:
\begin{eqnarray*}
&\phi (l,b,x,w,i): &bxx^*\rightarrow bxx^*\1\stackrel{(1,2,3)\diamond \Lambda _w}
   {\longrightarrow }(x^*bx)(w^*w)
   \stackrel{id\diamond (\epsilon _i\diamond id_w)}{\longrightarrow }\\
   &&\rightarrow (x^*bx)(x^*b^lw)\rightarrow x^*(b(xx^*)b^l)w
   \stackrel{\lambda _x}{\longrightarrow } x^*(bb^l)w, \\
&\psi (l,b,x,w,i): &x^*(b^lb)w\rightarrow (x^*b^l)bw 
   \stackrel{\epsilon _i^{*}}{\longrightarrow } w^*bw
   \stackrel{(1,3,2)}{\longrightarrow }bww^*.
\end{eqnarray*}

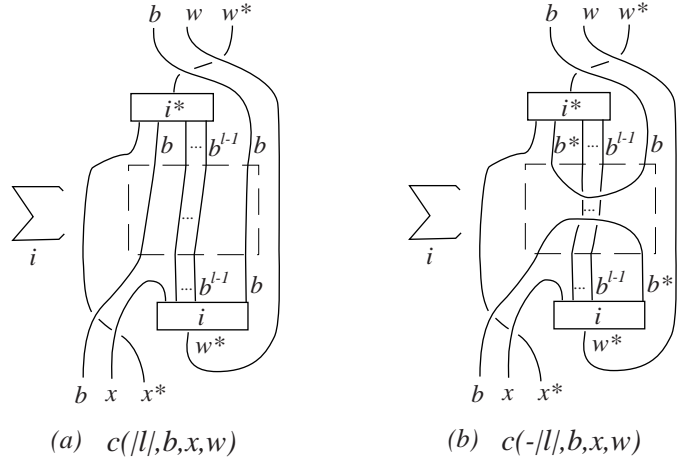
\begin{figure}[h]
\setlength{\unitlength}{1cm}
\begin{center}
\begin{picture}(9,6)
\put(0,0){\fig{figcirc.pict}}
\end{picture}
\end{center}
\caption{The morphism $c(l,b,x,w)$.}
\label{circ.pow}
\end{figure}

\begin{figure}[h]
\setlength{\unitlength}{1cm}
\begin{center}
\begin{picture}(15,8)
\put(0,-0.3){\fig{circind.pict}}
\end{picture}
\end{center}
\caption{Inductive step in the proof of proposition \ref{powcirc}.}
\label{ind.circ}
\end{figure}
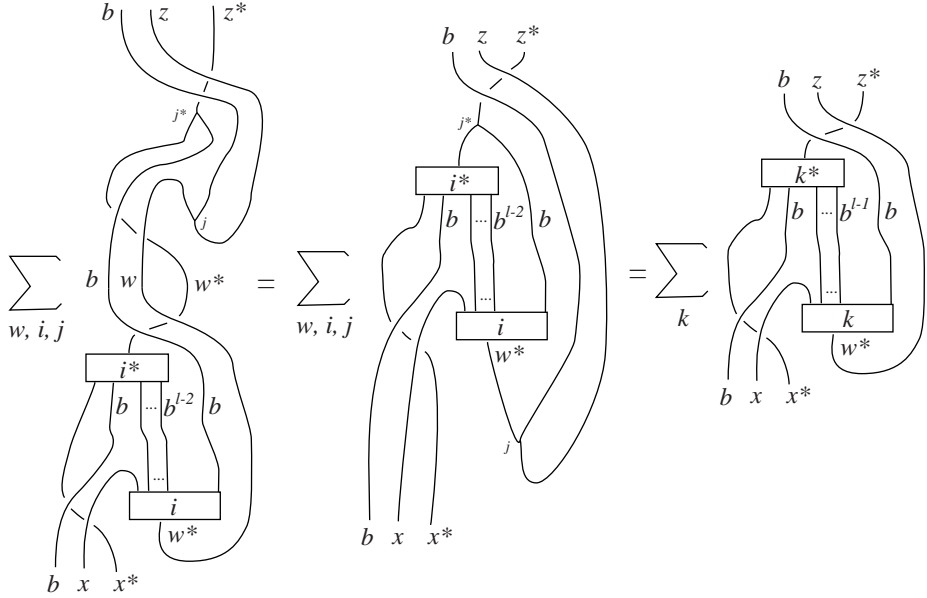

\begin{prop}
\label{powcirc} Given $f\in F(a,bxx^*)$, 
$
CR(b,a)^l(f)=\oplus _{w\in\Sigma} cr(l,b,x,w)\circ f,
$
where $cr(l,b,x,w)=\sum _i \psi(l,b,x,w,i)\circ (id_{x^*}\diamond 
f(b,l)\diamond id_w)\circ \phi (l,b,x,w,i)$.
\end{prop}
The corresponding diagram is shown in figure \ref{circ.pow}.
Note that $cr(l,b,x,w)$ can actually be written as the composition of
a partial trace and another morphism, and therefore it is independent 
on the particular choice of the basis $\{\epsilon _i\}_{i}$.
The fact that the circulator $CR(a,b):F(a,bxx^*)\rightarrow F(a,bww^*)$ 
acts as $cr(1,b,x,w)$
follows from the comparison of the diagrams in figure \ref{def.circ} and
\ref{circ.pow} (a). Then in figure \ref{inv.circ} we prove that 
$cr(-1,b,x,w)$ is a right inverse of the circulator. In a similar 
fashion one can see that it is a left inverse as well. 
Then, by assuming that for $l>1$, the $(l-1)$-th power of the circulator acts as
$cr(l-1,b,x,w)$, figure \ref{ind.circ} proves that the
statement is also true for the $l$-th power. 
The inductive proof for the negative 
powers of the circulator goes in a similar fashion. 
\begin{figure}[h]
\setlength{\unitlength}{1cm}
\begin{center}
\begin{picture}(9,11)
\put(0,-0.3){\fig{figcircinv1.pict}}
\end{picture}
\end{center}
\caption{Proof that $c(-1,b,x,w)$ represents the inverse of $c(1,b,x,w)$.}
\label{inv.circ}
\end{figure}
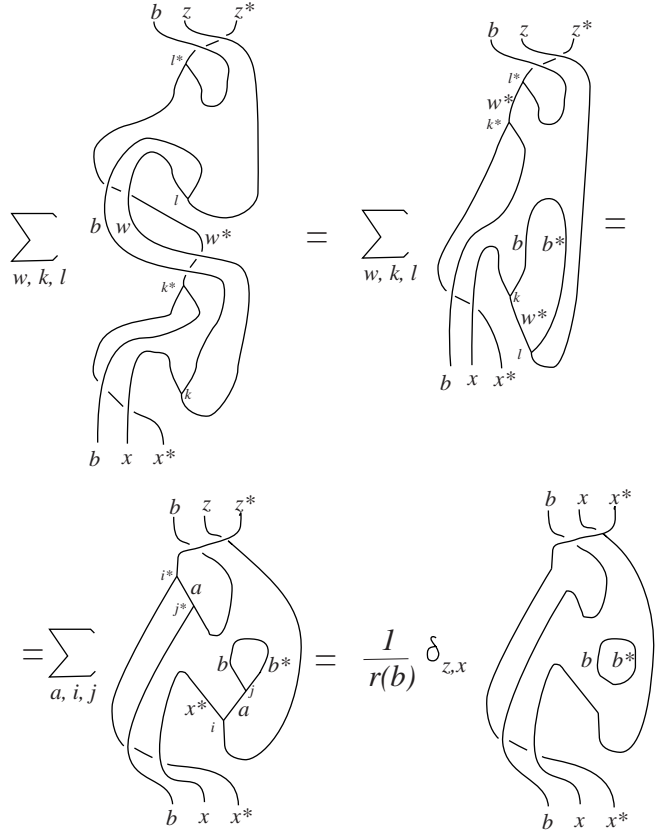 

\begin{figure}[h]
\setlength{\unitlength}{1cm}
\begin{center}
\begin{picture}(15,8)
\put(0,-0.3){\fig{figindmove.pict}}
\end{picture}
\end{center}
\caption{Inductive step in the proof of theorem \ref{qequiv}.}
\label{endind}
\end{figure}
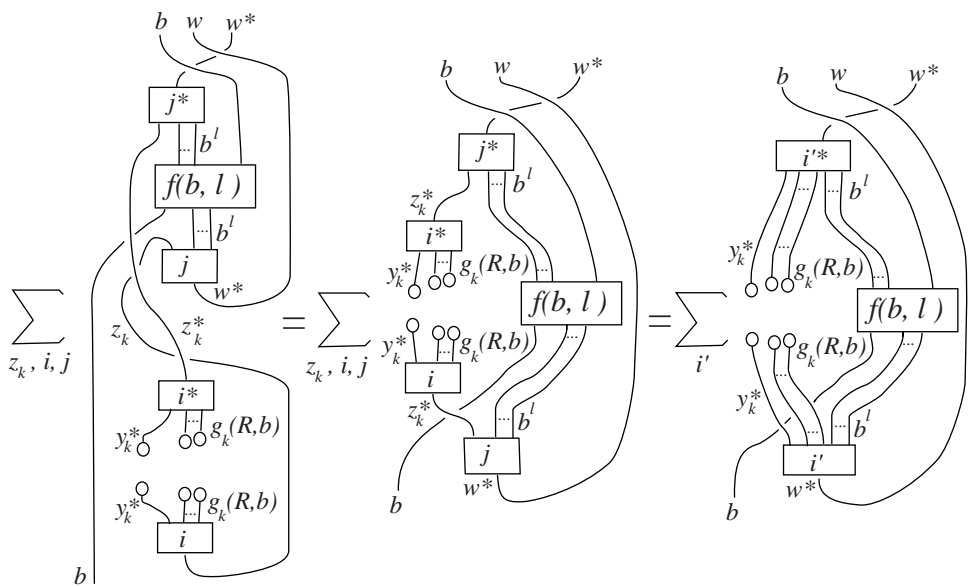

Now theorem \ref{qequiv} follows from the following statement.
Let $P$ be a presentation with $n$ generators,  
$R=x_{i_{1}}^{l_{1}}x_{i_{2}}^{l_{2}}\ldots x_{i_{s}}^{l_{s}}$ be
a relation in $P$, and 
given  $\underline{y},\underline{z}\in\Sigma ^{\times n}$ and
$b\in\Sigma$, let
 $\{\epsilon _{j}^{i}(z^*_{i},y^*_{i}g_i(R,b))\}_{j}$, $i=1\div n$,
  be a set of bases for
the homomorphisms spaces $F(z^*_{i},y^*_{i}g_i(R,b))$. Then the map 
$V(n+1)\rightarrow V(n+1)$, which the relation $R$ induces according to
the algorithm in \ref{algorithm}, is given by
$
\varphi \rightarrow \oplus _{\underline{y},\underline{z},j}
                rl(\underline{y},\underline{z},j)\circ \varphi,
$
where $rl(\underline{y},\underline{z},j)$ is the following composition
of maps:
\begin{eqnarray*}
&&b^{*}by_{1}y_{1}^{*}y_{2}y_{2}^{*}\dots y_{n}y_{n}^{*}
    \rightarrow 
   b^{*}(\diamond _{i}y_{i}^{*})b(\diamond _{i}y_{i})
   \stackrel{id\diamond (\diamond _{i}\Lambda 
   _{z_{i}})}{\longrightarrow}\\
&&b^{*}(\diamond _{i}y_{i}^{*})b(\diamond _{i}y_{i})
   (\diamond _{i}z_{i}^{*}z_{i})
   \stackrel{id\diamond (\diamond _{i}(\epsilon ^{i}_{j}\diamond 
   id_{z_{i}}))}
      {\longrightarrow}\\
&& b^{*}(\diamond _{i}y_{i}^{*})b(\diamond _{i}y_{i})
   (\diamond _{i}y_{i}^{*}g_i(R,b)z_{i})
   \rightarrow 
   b^{*}(\diamond _{i}y_{i}^{*})b(\diamond _{i}y_{i}
   y_{i}^{*})(\diamond _{i}g_i(R,b))(\diamond _{i}z_{i})
   \stackrel{id\diamond (\diamond _{i}\lambda _{y_{i}})\diamond id}
      {\longrightarrow}\\
&&b^{*}(\diamond _{i}y_{i}^{*})b(\diamond _{i}g_i(R,b))(\diamond _{i}z_{i})
   \stackrel{id\diamond \phi(R)\diamond id}{\longrightarrow}
  b^{*}(\diamond _{i}y_{i}^{*})(\diamond _{i}g_i(R,b))b(\diamond _{i}z_{i})
   \rightarrow\\
&&b^{*}(\diamond _{i}y_{i}^{*}g_i(R,b))b(\diamond _{i}z_{i})
   \stackrel{id\diamond (\diamond _{i}\epsilon ^{i}_{j}{}^*)\diamond id}
      {\longrightarrow}
  b^{*}(\diamond _{i}z^{*}_{i})b(\diamond _{i}z_{i})
   \rightarrow
   b^{*}bz_{1}z_{1}^{*}z_{2}z_{2}^{*}\dots z_{n}z_{n}^{*},
\end{eqnarray*}
where $\phi(R)=(id_{b}\diamond \kappa (R))[R,b]^{o}(\kappa 
      (R)^{-1}\diamond id_{b}):b(\diamond _{i}g_i(R,b))\rightarrow
      (\diamond _{i}g_i(R,b))b$ (see \ref{onerel}).  
The 
corresponding diagram is presented in figure \ref{endstep}, where
for simplicity we have restricted ourselves to the case of two generators.
The statement is proved by induction over $s$. For $s=1$ it is
reduced to proposition \ref{powcirc}. Assume it is true for $R$, and
let $R'=Rx_{k}^{l}$ for some $k$ and $l$. Then
\begin{eqnarray*}
&g_{i}(R',b)=\left\{\begin{array}{c}g_{k}(R,b)b^{l} \mbox{ if } i=k,\\
                                    g_i(R,b) \mbox{ otherwise }
                   \end{array}\right. , 
&\mbox{ and }\quad [R',b]=\circ _{b}([R,b]^{o}\circ _{b}f(b,l)).
\end{eqnarray*}
Then the statement follows from figure 
\ref{endind}.                  

\begin{figure}[h]
\setlength{\unitlength}{1cm}
\begin{center}
\begin{picture}(12,9)
\put(0,0){\fig{figendmove.pict}}
\end{picture}
\end{center}
\caption{The morphisms $V(n+1)\rightarrow V(n+1)$ induced by a 
          relation $R$ (case $n=2$).}
\label{endstep}
\end{figure}
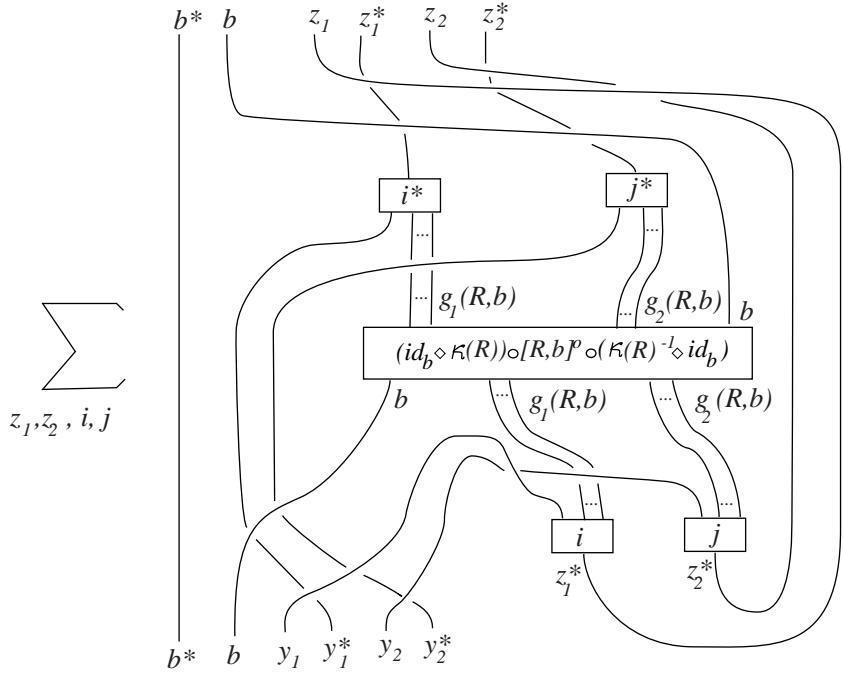 

\section{Conjectures and consequences}

As we said, this work was inspired by an extensive numerical study of the
invariants defined in \cite{Q:lectures}. The numerical project has been carried
out by Frank Quinn, the author, and Luoqi Zhang, and  the autonomous 
tensor categories 
used in it are subcategories of the Gelfand-Kazhdan 
categories, as described in \ref{perp.cat.}. We remind that
such a subcategory is determined by a type and rank of simple Lie algebra 
($A,B,C,D,E$ or $F$), a prime number $p$, and an invariant 
sublattice of the weight lattice (containing the root lattice). 
But actually many of the categories corresponding to different algebras
are equivalent. Conjecturally,
any two autonomous tensor categories with sets of simple objects 
$\Sigma _1$ and $\Sigma _2$, such that there is a bijection 
$\mu :\Sigma _1\rightarrow \Sigma _2$, and $dim(F(a,bc))=
dim(F(\mu(a),\mu(b)\mu(c)))$,
are equivalent. In other words, the dimensions of the homomorphism
spaces determine the category. Let $dim(a,bc)$ denote the
dimension of $F(a,bc)$ and we will refer to those as dimension functions. 
On web page http://www.math.vt.edu/quantum\_topology has been collected
a list with the values of the dimension functions for some categories
corresponding to algebras type $A,B,C,D$, small primes, and either the
full weight lattice or the root lattice. Let $L$ be the order of the
abelian group obtained as a quotient of the lattice over which the category is
defined, modulo the root lattice. Then all numerically generated 
examples satisfy the following conjecture:
\begin{conj} 
\label{conj1}In the case of stable Gelfand-Kazhdan categories, 
for any $b\in\Sigma$
\begin{itemize}
\item[(a)] 
$dim(a,b^{Lp})= 0 \quad (\mbox{mod $ p$}),\mbox{ if } a\neq \1 $;

\item[(b)] 
$r(b)^{-1}\sum _{c\in\Sigma}dim(b,cc^*) =\left\{
    \begin{array}{c} |\Sigma | \mbox{ if $b$ belongs to the root lattice,}\\
                     0 \mbox{ otherwise} 
    \end{array} \quad (\mbox{mod $ p$})\right.$.
\end{itemize}
\end{conj}
We make few comments on this. From \ref{rank} and (a) above we obtain
$$
dim(\1,b^{Lp})=\sum _{a\in\Sigma}dim(a,b^{Lp})r(a)=r(b^{Lp})=
r(b)^{Lp}=r(b)^L, \quad (mod\; p).
$$
In connection with (b) we comment that the statement can be proved 
easily for algebra type $A_1$. Moreover, the 0 part can be derived from a 
conjecture of Frank Quinn, according to which the category over
an arbitrary lattice is a product of the category over the root lattice
and a finite group category. If this is true, it would imply that the
only interesting categories are the ones coming from the root lattice,
since a finite group category brings to a classical invariant.
For categories over the root lattice, the conjecture \ref{conj1} brings
to the following result.

\begin{corr} 
\label{corr}Let $P=<x_1,x_2, \dots ,x_n|R_1,R_2,\dots ,R_m>$ be a
group presentation. If we are working with a Gelfand-Kazhdan category
over the root lattice,
and \ref{conj1} (b) is true, then
$
Q(P')=|\Sigma|\, Q(P''),
$
where 
\begin{eqnarray*}
&&P'=<x_1,x_2, \dots ,x_n,y|R_1,R_2,\dots ,R_m, x_kyx_k^{-1}y^{-1}>\mbox{ and}\\
&&P''=<x_1,x_2, \dots ,x_n|R_1,R_2,\dots ,R_m, x_k>.
\end{eqnarray*}
\end{corr}

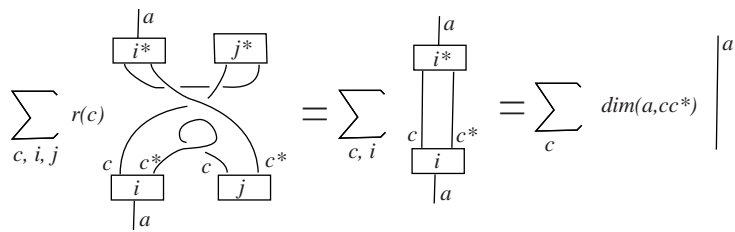
\begin{figure}[h]
\setlength{\unitlength}{1cm}
\begin{center}
\begin{picture}(10,3)
\put(0,-0.3){\fig{figcorr.pict}}
\end{picture}
\end{center}
\caption{Proof of corollary \ref{corr}.}
\label{fig.corr}
\end{figure}

Let
$\underline{b}\in\Sigma ^{\times m}$, and $a,c\in\Sigma$. 
If $(\underline{b},c)$ denotes the sequence $(b_{1},b_{2}, \ldots 
,b_{m},c)$ we have that
\begin{eqnarray*}
&&[P',(\underline{b},c)]=
[P,\underline{b}]\diamond 
   [x_kyx_k^{-1}y^{-1},c]:Rel(\underline{b})ccc^{*}c^{*}\rightarrow 
    Rel(\underline{b})ccc^{*}c^{*},\\
&&G_l'(\underline{b},c)=\left\{\begin{array}{l}
    G_k(\underline{b}) (cc^*)\mbox{ if }l=k,\\
     G_l(\underline{b})\mbox{ otherwise, 
     }\end{array}\right.\quad\mbox{ and }\quad
         G_y'(\underline{b},c)=cc^*, \\
&&\xi (P'): Rel(\underline{b})ccc^*c^*
  \stackrel{\xi (P)\diamond (2,3)}{\longrightarrow}
  (\diamond _{l}G_{l}(\underline{b}))cc^*cc^*
  \rightarrow \\
&&\qquad \qquad\qquad\rightarrow 
  (\diamond _{l<k}G_{l}(\underline{b}))G_k(\underline{b}) (cc^*)
  (\diamond _{l>k}G_{l}(\underline{b}))cc^*.
\end{eqnarray*}
Moreover,
\begin{eqnarray*}
&&[P'',(\underline{b},a)]= r(a)
[P,\underline{b}]\diamond 
   id_{a}:Rel(\underline{b})a\rightarrow 
    Rel(\underline{b})a,\\
&&G_l''(\underline{b},a)=\left\{\begin{array}{l}
    G_k(\underline{b}) a\mbox{ if }l=k,\\
     G_l(\underline{b})\mbox{ otherwise, }\end{array}\right.\\
&&\xi (P'')= Rel(\underline{b})a
  \stackrel{\xi (P)\diamond id_{a}}{\longrightarrow}
  (\diamond _{l}G_{l}(\underline{b}))a
  \rightarrow
  (\diamond _{l<k}G_{l}(\underline{b}))G_k(\underline{b}) a
  (\diamond _{l>k}G_{l}(\underline{b})),
\end{eqnarray*}
where we have used the notations in \ref{morerel}.
Then by using \ref{ptrprop} (c) and (d) we obtain the following for
the part of the partial trace 
involving $G_k'(\underline{b},c)$ and $G_y'(\underline{b},c)$:
\begin{eqnarray*}
&W&=\sum _{c\in\Sigma}{\tt Tr}_{\{G_k'(\underline{b},c)\rightarrow \1,
      \;G_y'(\underline{b},c)\rightarrow \1\}}
 (\xi (P')\circ ([P,\underline{b}]\diamond [x_kyx_k^{-1}y^{-1},c])\circ
\xi (P')^{-1})=\\
&&=\sum _{c,a\in\Sigma}{\tt Tr}_{G_k(\underline{b})a\rightarrow \1}(
\xi (P'')\circ ([P,\underline{b}]\diamond K)\circ\xi (P'')^{-1}),
\end{eqnarray*}
where $K={\tt Tr}_{\{cc^*\rightarrow a, cc^*\rightarrow \1\}}(
 (2,3)\circ[x_kyx_k^{-1}y^{-1},c]\circ (2,3))$.
Figure \ref{fig.corr} shows that in fact, $K=dim(a,cc^*)id_a\diamond id_{\1}$.
Hence we obtain :
\begin{eqnarray*}
&W&=\sum _{c,a\in\Sigma}dim(a,cc^*){\tt Tr}_{G_k(\underline{b})a\rightarrow \1}(
\xi (P'')\circ ([P,\underline{b}]\diamond id_{a})\circ\xi (P'')^{-1})=\\
&&=|\Sigma |\sum _{a\in\Sigma}{\tt Tr}_{G_k(\underline{b})a\rightarrow \1}(
\xi (P'')\circ [P'',(\underline{b},a)]\circ\xi (P'')^{-1}),
\end{eqnarray*}
which is exactly the part of the partial trace in the definition of 
$Q(P'')$ involving $G_k''(\underline{b},a)$. 

Another conjecture we would like to discuss concerns the order of the
circulator.  It is observed in \cite{Q:lectures}, that for all numerically
generated categories, the order of the circulator is in fact $Lp$. 
On the basis of this and the expression for the powers of the circulator
\ref{powcirc} we may state the following 

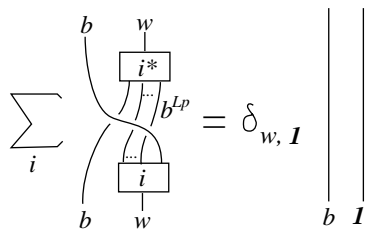
\begin{figure}[h]
\setlength{\unitlength}{1cm}
\begin{center}
\begin{picture}(4,3)
\put(0,0){\fig{figconj.pict}}
\end{picture}
\end{center}
\caption{Conjecture \ref{conj2}.}
\label{fig.conj2}
\end{figure}

\begin{conj} 
\label{conj2}
Let $\sigma^{Lp+1}(b)=(1,2,3,\dots ,Lp+1):bb^{Lp}\rightarrow
bb^{Lp}$. Then in the case of  stable  Gelfand-Kazhdan
categories,
$$
{\tt Tr}_{b^{Lp}\rightarrow w}\sigma^{Lp+1}(b)=\left\{
\begin{array}{ll}
   0 &\mbox{if }\quad w\neq \1,\\
   id_b \diamond id _{\1} &\mbox{otherwise.}
\end{array} \right.
$$
\end{conj}
The corresponding diagram is presented in figure \ref{fig.conj2}.

The conjecture about the finite order of the circulator implies that given
a presentation, the exponents of the generators in the relations matter only
mod $Lp$, in particular it is enough to look at presentations $P$ with only
positive exponents. In this last case, the morphism $[P,\underline{b}]$
has a particularly simple form - it is a product of cyclic permutations:
$$
[P,\underline{b}]=(1,2,3,\dots ,e_1)\diamond (1,2,3,\dots ,e_2)
\diamond \dots \diamond (1,2,3,\dots ,e_m),
$$
where $e_i$ is the total exponent of the relation $R_i$. This fact doesn't 
seem to help though with evaluating the invariant, at least until some
more information is available for the symmetry group representations
$\rho [b^k]$. 

In conclusion we point out that clearly the main question is how the
invariant behaves under the move (vii) in section 1, i.e. under simple
homotopy equivalence. The numerical examples haven't exhibited any
nontrivial behavior, but using the formalism in this paper to 
write down the change in the value of the
invariant under (vii), brings to
expressions which look quite complicated and we don't know how to analyze. 
We suspect that if the behavior of the computed invariant under this 
move is 
actually trivial, this is a consequence of some specific properties of the 
Gelfand - Kazhdan
categories being used ( like the ones listed in conjecture \ref{conj1}), 
rather then a consequence of the general setup of the theory.

\bibliographystyle{amsalpha}

\end{document}